\documentclass[twocolumn]{autart}    

\newcommand{\R}{{\mathbb R}}  
\newcommand{\Rn}{\R^{n}}

\usepackage{graphicx}
\usepackage{amsmath}
\usepackage{amsfonts}
\usepackage{subcaption}
\usepackage{amssymb}
\usepackage{url} 

\newcommand\norm[1]{\lVert#1\rVert}
\newcommand\Tr[1]{\mathrm{Tr}\left(#1\right)}
\newcommand\grad[1]{\mathrm{grad}\left(#1\right)}
\newcommand\eigmax[1]{\lambda_{\mathrm{max}}\left(#1\right)}
\newcommand\eigmin[1]{\lambda_{\mathrm{min}}\left(#1\right)}
\newcommand{\de}{\mathrm{d}}
\newcommand\innprod[2]{\langle #1, \, #2 \rangle}

\usepackage{xcolor} 
\newcommand{\Review}[1]{\textcolor{black}{#1}}

\makeatletter
\renewenvironment{pf}%
  {\par\addvspace{\@bls \@plus 0.5\@bls \@minus 0.1\@bls}\noindent
   {\bfseries\Elproofname}\enspace\ignorespaces}%
  {\enspace\hfill$\square$\par\addvspace{\@bls \@plus 0.5\@bls \@minus 0.1\@bls}} 
\makeatother

\begin{document}

\begin{frontmatter}

\title{Perturbed Gradient Descent Algorithms are \\ Small-Disturbance Input-to-State Stable} 

\thanks[footnoteinfo]{This work was supported by the U.S. Food and Drug Administration under the FDA BAA-22-00123 program (Award Number 75F40122C00200) for L.C. and R.D.B., the National Science Foundation (NSF) grants CNS-2227153 and ECCS-2210320 for Z.P.J., and the Air Force Office of Scientific Research (AFOSR) grant FA9550-21-1-0289 and the Office of Naval Research (ONR) grant N00014-21-1-2431 for E.D.S. The material in this paper was not presented at any conference.}

\author[LC]{Leilei Cui}\ead{lcui@unm.edu},    
\author[ZPJ]{Zhong-Ping Jiang}\ead{zjiang@nyu.edu},               
\author[EDS]{Eduardo D. Sontag}\ead{e.sontag@northeastern.edu},   and   
\author[RB]{Richard D. Braatz}\ead{braatz@mit.edu}  

\address[LC]{Department of Mechanical Engineering, University of New Mexico, Albuquerque, NM, USA}  
\address[ZPJ]{Department of Electrical and Computer Engineering and Department of Civil and Urban Engineering, New York
University, Brooklyn, NY, USA}        
\address[EDS]{Department of Electrical and Computer Engineering and Department of BioEngineering, Northeastern University,
Boston, MA, USA}             
\address[RB]{Massachusetts Institute of Technology, Cambridge, MA, USA} 
          

\begin{abstract}                          
This article investigates the robustness of gradient descent algorithms under perturbations. The concept of small-disturbance input-to-state stability (ISS) for discrete-time nonlinear dynamical systems is introduced, along with its Lyapunov characterization. The conventional linear \emph{Polyak-\L{}ojasiewicz} (PL) condition is then extended to a nonlinear version, and it is shown that the gradient descent algorithm is small-disturbance ISS provided the objective function satisfies the generalized nonlinear PL condition. This small-disturbance ISS property guarantees that the gradient descent algorithm converges to a small neighborhood of the optimum under sufficiently small perturbations. As a direct application of the developed framework, we demonstrate that the LQR cost satisfies the generalized nonlinear PL condition, thereby establishing that the policy gradient algorithm for LQR is small-disturbance ISS. Additionally, other popular policy gradient algorithms, including natural policy gradient and Gauss-Newton method, are also proven to be small-disturbance ISS.
\end{abstract}

\end{frontmatter}

\section{Introduction}
\Review{
Gradient-based optimization algorithms are central to modern machine learning, enabling efficient search over high-dimensional parameter spaces. Methods such as gradient descent and adaptive moment estimation (Adam) \cite{kingma2014adam} are among the most widely used first-order optimizers in deep learning. Understanding the convergence behavior of gradient descent is therefore essential. While classical analyses assume exact gradient information, in practice gradients are often perturbed by round-off errors, noisy measurements, inaccurate formulas, or approximations in auxiliary computations (see \cite[Chapter 4]{book_Polyak} and \cite[p.\ 38]{bertsekas1997nonlinear}). These perturbations may cause oscillations near the optimum or even divergence \cite{book_Polyak,bertsekas1997nonlinear}. Consequently, beyond guaranteeing convergence in noise-free settings, gradient-based optimization algorithms should be designed to maintain robustness under perturbations, motivating the joint study of convergence and robustness properties.
}

\Review{
A useful perspective for understanding optimization is to interpret gradient-based algorithms as dynamical systems. This viewpoint enables the application of control-theoretic tools, such as Lyapunov stability, to analyze algorithm behavior \cite{Ashia2021,Bin2017,moucer2023systematic,jordan2018dynamical}. However, Lyapunov stability primarily characterizes system behavior in the absence of external inputs and is therefore limited in assessing robustness under perturbations. Input-to-state stability (ISS) extends Lyapunov stability by explicitly relating the system state to the magnitude of external inputs, providing a more suitable framework for analyzing the impact of perturbations on gradient-based algorithms \cite{Sontag1989,Sontag2008,Sontag2022,CJS2024}. ISS characterizes both asymptotic stability in noise-free settings and the size of the neighborhood reached under persistent disturbances, making it a valuable tool for studying both convergence and robustness of gradient-based methods. For example, \cite{ieee_tac_2018_cherukuri_et_al_convexity_saddle_point_dynamics} shows that gradient descent–ascent dynamics for convex–concave functions are ISS with respect to additive noise. In \cite{kolmanovsky2022inputtostate}, ISS is used to analyze the robustness of a bilevel optimization algorithm subject to errors arising from incomplete inner-loop computations. ISS has also been applied to extremum-seeking methods \cite{2021_poveda_krstic_fixedtime_iss_extremum_seeking,2021_arxiv_iss_gradient_suttner_dashkovskiy}, as well as to output regulation problems involving tracking of gradient flows in systems subject to plant-level disturbances \cite{2020arxiv_bianchin_poveda_dallanese_gradient_iss_switched_systems}.
}

\Review{
In this paper, we establish a connection between the ISS property of gradient descent algorithms and the \emph{Polyak-\L{}ojasiewicz (PL)} condition \cite{Polyak1963,Lojasiewicz1963,Hamed2016}. The PL condition guarantees linear convergence of gradient descent without requiring convexity by imposing $\norm{\nabla \mathcal{J}(z)}^2 \ge c(\mathcal{J}(z)-\mathcal{J}^*)$ for some constant $c>0$, where $\mathcal{J}(z)$ is continuously differentiable on $\mathcal{Z}$ with minimum value $\mathcal{J}^*$. However, the linearity imposed by the conventional PL condition restricts its applicability to a broader range of problems, such as the optimization of LQR costs discussed in the second part of the paper, where globally establishing such linearity is infeasible. This limitation motivates generalizing the PL condition to nonlinear forms that broaden applicability while preserving convergence and robustness guarantees for gradient descent algorithms.
}

\Review{
A natural generalization replaces the constant $c$ in the PL condition with a function $\alpha$ belonging to class $\mathcal{K}$, which requires $\alpha$ to be strictly increasing and vanish at zero. The resulting condition $\norm{\nabla \mathcal{J}(z)}^2 \ge \alpha(\mathcal{J}(z)-\mathcal{J}(z^*))$ is referred to as a $\mathcal{K}$-PL condition. For example, $\alpha(r)=\frac{r}{1+r}$ is a valid $\mathcal{K}$-function that saturates as $r\to\infty$. We show that if the objective function is coercive, has an $L$-Lipschitz continuous gradient, and satisfies the $\mathcal{K}$-PL condition, then gradient descent with step size $0<\eta\le 1/L$ is \emph{small-disturbance ISS}. In this case, the gradient descent converges to a neighborhood of the optimum whose size scales with the magnitude of the perturbation. If the $\mathcal{K}$-PL condition is further strengthened to a $\mathcal{K}_\infty$-PL condition---requiring $\alpha(r) \to \infty$ as $r \to \infty$---the gradient descent algorithm is ISS \cite{Sontag2022}. If the $\mathcal{K}$-PL condition is relaxed to a $\mathcal{PD}$-PL condition, which requires $\alpha$ to be positive definite, the gradient descent algorithm is integral ISS, which is equivalent to global asymptotic stability for discrete-time dynamical systems \cite{angeli1999intrinsic}. Meanwhile, the conventional (linear) PL condition implies that the gradient descent algorithm is exponentially ISS. It should be noted that the robustness of gradient descent algorithms is typically analyzed under the assumption of convexity \cite{book_Polyak,solodov1998error,devolder2014first,lessard2016analysis,Cyrus2018,Aybat2020robust,Ashia2021,Mihailo2024} or by assuming that the perturbations converge to zero \cite{Bertsekas2000}. In contrast, the generalized $\mathcal{K}$-PL condition offers an alternative framework for robustness analysis when these convexity or convergent-perturbation assumptions do not hold.
}

\Review{
A direct application of the proposed $\mathcal{K}$-PL condition is the robustness analysis of reinforcement learning (RL) algorithms for the linear quadratic regulator (LQR). Policy optimization provides an effective framework for RL by parameterizing policies with universal approximators and updating parameters via gradient descent \cite[Chapter 13]{book_sutton}. Introduced by Kalman \cite{kalman1960contributions}, LQR is a fundamental control method widely used in engineering. When system matrices are known, both the gradient and the optimal solution can be computed explicitly, allowing performance comparisons with the optimal policy. Consequently, LQR has long served as a benchmark problem in control theory \cite{1970levine_athans,makila1987computational} and has recently gained renewed attention in RL \cite{fazel2018global,Mohammadi2022,Hu2023Review}. In model-free settings, gradients must be estimated from data, for example using finite-difference methods \cite{fazel2018global,Mohammadi2022,Li2021} or adaptive dynamic programming approaches \cite{books_Bertsekas,book_Jiang}. Such estimation inevitably introduces errors due to noisy measurements and limited samples. Therefore, beyond convergence in noise-free scenarios, robustness with respect to gradient errors is essential for understanding RL algorithms. Using coercivity, Lipschitz continuity of the gradient, and the $\mathcal{K}$-PL property of the LQR cost \cite{CJS2024}, we show that standard gradient descent for LQR is small-disturbance ISS. We further establish that natural gradient descent and Gauss-Newton gradient descent algorithms for LQR also exhibit small-disturbance ISS.
}

\Review{
While the concept of small-disturbance ISS for discrete-time dynamical systems stems from the continuous-time framework \cite{CJS2024}, its Lyapunov characterization cannot be directly established by simply discretizing the arguments used in the continuous case. Since gradient descent algorithms inherently operate in discrete time, it is necessary to present these results with comprehensive and rigorous proofs. Moreover, analyzing the convergence of discrete-time gradient descent is more challenging due to its dependence on the Lipschitz constant of the gradient and the careful design of the step size. For the standard policy gradient descent for LQR, small-disturbance ISS is ensured by carefully selecting the step size based on the Lipschitz continuity of the gradient. For natural gradient descent and Gauss-Newton gradient descent algorithms, the step sizes are designed by directly differencing the two Lyapunov functions corresponding to consecutive updates, thereby ensuring small-disturbance ISS.
}

\Review{
In summary, the contributions of this paper are threefold. First, we propose a Lyapunov-like necessary and sufficient condition for small-disturbance ISS in discrete-time dynamical systems. Second, under coercivity, the $\mathcal{K}$-PL condition, and $L$-Lipschitz continuity of the objective function and its gradient, we show that gradient descent with step size $0 < \eta \le \frac{1}{L}$ is small-disturbance ISS. Third, by analyzing the LQR cost, we establish that the standard policy gradient algorithm, as well as the natural gradient and Gauss-Newton gradient methods, are all small-disturbance ISS.
}

\Review{
The remainder of the paper is organized as follows. Section 2 introduces notation and preliminaries. Section 3 presents the definition and Lyapunov-like characterization of small-disturbance ISS. Section 4 establishes small-disturbance ISS of the perturbed gradient descent algorithm. Section 5 shows that the standard, natural, and Gauss–Newton gradient descent methods for LQR are all small-disturbance ISS. Section 6 presents simulation results. Section 7 concludes the paper.
}
\section{Notations and Facts}

In this article, we denote by $\mathbb{R}$ ($\mathbb{R}_+$) the set of (nonnegative) real numbers, $\mathbb{Z}_+$ the set of nonnegative integers, and $\mathbb{S}^n$ ($\mathbb{S}^n_{+}/\mathbb{S}^n_{++}$) the set of $n$-dimensional real symmetric (positive semidefinite/definite) matrices. For a real symmetric matrix, $\eigmin{\cdot}$ and $\eigmax{\cdot}$ denote the minimum and maximum eigenvalues, respectively. $\Tr{\cdot}$ is the trace of a square matrix. The Euclidean norm of a vector or spectral norm of a matrix is denoted by $\norm{\cdot}$, while $\norm{\cdot}_F$ denotes the Frobenius norm of a matrix.

Let $\ell_\infty^{n}$ ($\ell_\infty^{m \times n}$) denote the set of bounded functions $w: \mathbb{Z}_+ \to \Rn$ ($K: \mathbb{Z}_+ \to \mathbb{R}^{m \times n}$), with the $\ell_\infty$-norm given by $\norm{w}_\infty = \sup_{k \in \mathbb{Z}_+}\norm{w(k)}$ ($\norm{K}_\infty = \sup_{k \in \mathbb{Z}_+}\norm{K(k)}_F$). The truncation of $w \in \ell_\infty^{n}$ at step $k$ is denoted by $w_{[k]}$, defined as
$$
w_{[k]}(j) = \begin{cases}
      w(j) & \text{if } j \le k\\
      0 &  \text{if } j > k
    \end{cases}
$$
We use $I_n$ to denote the $n$-dimensional identity matrix and $\mathrm{Id}$ for the identity function. For any $K_1, K_2 \in \mathbb{R}^{m \times n}$ and $Y \in \mathbb{S}^n_{++}$, we define the inner product as $\innprod{K_1}{K_2}_Y = \Tr{K_1  Y  K_2^\top}$. For simplicity, we write $\innprod{K_1}{K_2} = \innprod{K_1}{K_2}_{I_n}$. Note that for any $K \in \mathbb{R}^{m \times n}$, $\norm{K}_F^2 = \innprod{K}{K}$. Finally, for any $A \in \mathbb{S}^n$ and $B \in \mathbb{S}^n$, $A \succ B$ indicates that $A - B \in \mathbb{S}^n_{++}$, and $A \succeq B$ indicates that $A - B\in \mathbb{S}^n_{+}$.

The concepts of comparison functions \cite{book_Hahn}, which are essential for stability analysis, are introduced here. A function $\alpha: \mathbb{R}_+ \to \mathbb{R}_+$ is defined as a $\mathcal{K}$-function if it is continuous, strictly increasing, and equals zero at the origin. For any $d >0$, a function $\alpha: [0,d) \to \mathbb{R}_+$ is a $\mathcal{K}_{[0,d)}$-function if it is continuous, strictly increasing, and vanishes at zero. A function $\alpha: \mathbb{R}_+ \to \mathbb{R}_+$ is a $\mathcal{K}_\infty$-function if it meets the criteria of a $\mathcal{K}$-function and additionally satisfies $\alpha(r) \to \infty$ as $r \to \infty$. A function $\beta: \mathbb{R}_+ \times \mathbb{R}_+ \to \mathbb{R}_+$ is called a $\mathcal{KL}$-function if, for each fixed $t \geq 0$, $\beta(\cdot, t)$ is a $\mathcal{K}$-function, and for each fixed $r \geq 0$, $\beta(r, \cdot)$ is decreasing and approaches zero as $t \to \infty$.

Several established facts are introduced next to support the development of the main results in this paper.

\begin{lem}[Weak triangle inequality \cite{Jiang1994}]\label{lm:weakTriangle}
    For any $\mathcal{K}$-function $\alpha$, any $\mathcal{K}_\infty$-function $\rho$, and any $a,b \in \R_+$, it holds that
    \begin{align*}
        \alpha(a+b) \le \alpha \circ (\mathrm{Id} + \rho)(a) + \alpha \circ (\mathrm{Id} + \rho^{-1}) (b).
    \end{align*}
\end{lem}

\begin{lem}[Cyclic property of trace \cite{book_Petersen}]\label{lm:traceCyclic}
    For any $X,Y,Z \in \mathbb{R}^{n \times n}$, $\Tr{XYZ} = \Tr{ZXY} = \Tr{YZX}$.
\end{lem}

\begin{lem}[Trace inequality \cite{Wang1986}]\label{lm:traceIneq}
    For any $S\in \mathbb{S}^{n}$ and $P \in  \mathbb{S}^{n}_{+}$, it holds that
    \begin{align*}
        \eigmin{S} \Tr{P} \le \Tr{SP} \le \eigmax{S} \Tr{P}.
    \end{align*}
\end{lem}

\begin{lem}[Cauchy-Schwarz inequality]\label{lm:CSInequality}
    For any $K_1,K_2 \in \mathbb{R}^{m\times n}$, $R \in \mathbb{S}^{m}_{++}$, and $Y \in \mathbb{S}^n_{++}$, it holds that
    \begin{align}\label{eq:CSInequality}
        \innprod{K_1}{RK_2}_Y \leq \innprod{K_1}{RK_1}_Y^{1/2} \innprod{K_2}{RK_2}_Y^{1/2}. 
    \end{align}
\end{lem}
\begin{pf}
    The expression $\langle K_1, R K_2 \rangle_{Y}$ defines a valid inner product with respect to $K_1$ and $K_2$, from which the Cauchy–Schwarz inequality follows directly.
\end{pf}

\begin{lem}[Lemma 2.4 in \cite{CJS2024}]\label{lm:homomorphism}
   The map $h(v) = \frac{1}{1+\norm{v}}v: \mathbb{R}^m \to \{w\in \mathbb{R}^m| \norm{w} < 1\}$, with $h^{-1}(w) = \frac{1}{1-\norm{w}}w$, is a homeomorphism.
\end{lem}

\begin{lem}[Proposition 2.6 in \cite{Sontag2022}]\label{lm:radiallyUnbounded}
    Suppose $\omega_1, \, \omega_2: \mathbb{R}^n \to \mathbb{R}$ are continuous, positive definite with respect to $\chi^*$, and radially unbounded. Then, there exist $\mathcal{K}_\infty$-functions $\rho_1$ and $\rho_2$ such that
    \begin{align*}
        \rho_1(\omega_2(\chi)) \le \omega_1(\chi) \le \rho_2(\omega_2(\chi)), \quad \forall \chi \in \mathbb{R}^n.
    \end{align*}
\end{lem}

\begin{lem}[Theorem 18 in \cite{book_sontag}]\label{lm:LyaEquaIntegral}
    If $A \in \mathbb{R}^{n \times n}$ is Hurwitz, then the Lyapunov equation 
    \begin{align*}
        A^\top  P + PA + Q = 0
    \end{align*}
    has a unique solution for any $Q \in \mathbb{R}^{n \times n}$, and the solution can be expressed as
    \begin{align*}
        P = \int_{0}^\infty e^{A^\top  t} Q e^{At} \mathrm{d}t.
    \end{align*}
\end{lem}
The following corollary is a direct consequence of Lemma \ref{lm:LyaEquaIntegral} and the cyclic property of trace in Lemma \ref{lm:traceCyclic}.
\begin{cor}\label{cor:traceLya}
    Suppose that $A \in \mathbb{R}^{n \times n}$ is Hurwitz, $M,N\in \mathbb{S}^n$, and $P,Y \in \mathbb{S}^n$ are the solutions of
    \begin{align}
    \begin{split}
        &A^\top P + PA + M = 0 \\
        &AY + YA^\top + N = 0.
    \end{split}
    \end{align}
    Then, $\Tr{MY} = \Tr{NP}$.
\end{cor}


\section{Small-Disturbance Input-to-State Stability}
In this section, we investigate the dependence of state trajectories on the magnitude of the disturbances for the discrete-time nonlinear system:
\begin{align}\label{eq:nonLinearSys}
    \chi({k+1}) = f(\chi(k), w(k))
\end{align}
where $\chi(k) \in \mathcal{S}$ denotes the state evolving in an open subset $\mathcal{S} \subset \mathbb{R}^n$ which is homeomorphic to $\Rn$, $w \in \ell_\infty^{m}$ denotes the disturbance, and $f: \mathcal{S} \times \mathbb{R}^m \to \mathcal{S}$ is a continuous function. Assume that $\chi^* \in \mathcal{S}$ is the equilibrium of the unforced system, that is $\chi^* = f(\chi^*,0)$. Denote by $\chi(\cdot,\xi,w)$ the trajectory of system \eqref{eq:nonLinearSys} with the initial state $\chi(0) = \xi$ and disturbance $w \in \ell_\infty^{m}$.

Since system \eqref{eq:nonLinearSys} is defined in an open subset $\mathcal{S}$, instead of $\mathbb{R}^n$, a {\em size function} is introduced to assist in stability analysis and serves as a barrier function preventing escape from $\mathcal{S}$. 

\begin{defn}[\cite{Sontag2022}]\label{def:sizeFunc}
    A function $\mathcal{V}: \mathcal{S} \to \mathbb{R}$ is a {\em size function} for $(\mathcal{S},\chi^*)$ if $\mathcal{V}$ is 
    \begin{enumerate}
        \item continuous;
        \item positive definite with respect to $\chi^*$, i.e. $\mathcal{V}(\chi^*)=0$ and $\mathcal{V}(\chi)>0$ for all $\chi \neq \chi^*$,  $\chi \in \mathcal{S}$;
        \item coercive, i.e. for any sequence $\{\chi_k\}_{k=0}^\infty$, $\chi_k \to \partial \mathcal{S}$ or $\norm{\chi_k} \to \infty$, it holds that $\mathcal{V}(\chi_k) \to \infty$, as $k \to \infty$.
    \end{enumerate}
\end{defn}
For $\mathcal{S}=\mathbb{R}^n$, $\mathcal{V}(\xi) = \norm{\xi - \chi^*}$ is a natural choice of a size function. By resorting to size functions, we first introduce the concept of classic ISS.
\Review{
\begin{defn}[\cite{Sontag1989,Sontag2022}]
    The nonlinear system \eqref{eq:nonLinearSys} is ISS if there exist a size function $\mathcal{V}$, a $\mathcal{KL}$-function $\beta$, and a $\mathcal{K}_\infty$-function $\gamma$, such that for all inputs $w$, and all initial states $\chi(0) \in \mathcal{S}$, $\chi(k) $ remains in $\mathcal{S}$ and satisfies
    \begin{align}\label{eq:ISS}
         \mathcal{V}(\chi(k)) \leq \beta(\mathcal{V}(\chi(0)),k) + \gamma(\norm{w}_\infty), \quad \forall k \in \mathbb{Z}_+.
    \end{align}
\end{defn}
}
\Review{
ISS guarantees that the system state remains bounded for any bounded disturbance input, regardless of its magnitude. However, many systems (e.g., the policy gradient dynamics for LQR studied in the second half of this paper) exhibit ISS behavior only when the disturbance is sufficiently small, and may become unstable when the disturbance magnitude is large. Consequently, the classical ISS requirement can be overly restrictive for many practical systems. This observation motivates the development of a weaker notion, termed small-disturbance ISS, which characterizes robustness with respect to disturbances whose magnitude lies within a certain admissible range.
}
\begin{defn}[\cite{Pang_Jiang_2021,CJS2024}]\label{def:smallISS}
The nonlinear system \eqref{eq:nonLinearSys} is {\em small-disturbance input-to-state stable (ISS)} if there exist a size function $\mathcal{V}$, a constant $d > 0$ (possibly $\infty$), a $\mathcal{KL}$-function $\beta$, and a $\mathcal{K}_{[0,d)}$-function $\gamma$, such that for all inputs $w$ bounded by $d$ (i.e. $\norm{w}_\infty<d$), and all initial states $\chi(0) \in \mathcal{S}$, $\chi(k) $ remains in $\mathcal{S}$ and satisfies \eqref{eq:ISS}.
\end{defn}
By causality, the same definition would result for every $k \in \mathbb{Z}_+ \setminus \{0\}$ if $\norm{w}_\infty$ was replaced by $\norm{w_{[k-1]}}_\infty$ in \eqref{eq:ISS}, where $w_{[k-1]}$ is the truncation of $w$ at $k-1$. Additionally, the classical ISS definition can be recovered by setting $d = \infty$ and $\mathcal{V}(\chi(k)) = \norm{\chi(k)}$. Thus, small-disturbance ISS serves as an extension of the classical ISS \cite{Sontag1989,Jiang2001ISS}.

\begin{defn}\label{def:smallISSLyap}
A function $\mathcal{V}:\mathcal{S} \to \mathbb{R}$ is a {\em small-disturbance ISS-Lyapunov function} for system \eqref{eq:nonLinearSys} if 
\begin{enumerate}
    \item $\mathcal{V}$ is a size function for $(\mathcal{S},\chi^*)$;
    \item there exist a $\mathcal{K}$-function $\alpha_1$ and a $\mathcal{K}_\infty$-function $\alpha_2$ such that
    \begin{align}\label{eq:compaLya}
        \mathcal{V}(f(\xi,\mu)) - \mathcal{V}(\xi)\leq - \alpha_2(\mathcal{V}(\xi))
    \end{align}
    for any $\xi \in \mathcal{S}$ and $\mu \in \R^m$ so that $\norm{\mu} \le \alpha_1(\mathcal{V}(\xi))$.
\end{enumerate}
\end{defn}

\begin{rem}
    An equivalent property holds if the function $\alpha_2$ in \eqref{eq:compaLya} is only required to be continuous and positive definite \cite{Jiang2001ISS,Jiang2001ConverseLya}.
\end{rem}
The following remark provides a ``dissipation" type of characterization for the small-disturbance ISS property. 
\begin{rem}
    A size function $\mathcal{V}$ for $(\mathcal{S}, \chi^*)$ is a small-disturbance ISS-Lyapunov function for system \eqref{eq:nonLinearSys} if and only if there exist a $\mathcal{K}_\infty$-function $\alpha_2$, some $d>0$ (possibly $\infty$), and a $\mathcal{K}_{[0,d)}$-function $\alpha_3$ such that 
    \begin{align}\label{eq:smallISSLyapu_diss}
        \mathcal{V}(f(\xi,\mu)) - \mathcal{V}(\xi) \le - \alpha_2(\mathcal{V}(\xi)) + \alpha_3(\norm{\mu})
    \end{align}
    for all $\mu \in \R^m$ bounded by $d$, i.e., $\norm{\mu}<d$.
\end{rem}
Without loss of generality, assume that $\alpha_3(r) \to \infty$ as $r \to d$ (add $\frac{r}{d-r}$ to it if it is not). Clearly, \eqref{eq:smallISSLyapu_diss} implies \eqref{eq:compaLya}. Suppose now that Property 2 in Definition \ref{def:smallISSLyap} holds with some $\alpha_1 \in \mathcal{K}$ and $\alpha_2 \in \mathcal{K}_\infty$. Let $d = \sup_{r \in \R_+}\alpha_1(r)$ (possibly $\infty$). For any $r \in [0,d)$, define ${\alpha}_3(r) = \max\{\mathcal{V}(f(\xi,\mu)) - \mathcal{V}(\xi) + \alpha_2(\mathcal{V}(\xi))| \norm{\mu} \le r, \alpha_1(\mathcal{V}(\xi)) \le r\}$. Then, $\alpha_3$ is continuous, non-decreasing, and zero at zero. In addition, we can assume that $\alpha_3$ is a $\mathcal{K}_{[0,d)}$-function (add $\frac{r}{d-r}$ to it if it is not). When $\norm{\mu} \le \alpha_1(\mathcal{V}(\xi))$, $\mathcal{V}(f(\xi,\mu)) - \mathcal{V}(\xi) \le -\alpha_2(\mathcal{V(\xi)})$; when $\norm{\mu} \ge \alpha_1(\mathcal{V}(\xi))$, $\mathcal{V}(f(\xi,\mu)) - \mathcal{V}(\xi) \le -\alpha_2(\mathcal{V}(\xi)) + \alpha_3(\norm{\mu})$. Hence, \eqref{eq:smallISSLyapu_diss} holds.

As in classic Lyapunov stability theory, we can show that small-disturbance ISS is equivalent to the existence of a small-disturbance ISS-Lyapunov function.
\begin{thm}\label{thm:smallISSSufficient}
System \eqref{eq:nonLinearSys} is small-disturbance ISS if and only if it admits a small-disturbance ISS-Lyapunov function.
\end{thm}
\begin{pf} 
\textbf{Sufficiency:} The proof for the ISS property in \cite{Jiang2001ISS} is adapted here. We denote by $\chi(k) = \chi(k,\xi,w)$ the state trajectory of \eqref{eq:nonLinearSys} for fixed initial state $\xi \in \mathcal{S}$ and input $w \in \ell_\infty^m$. For the input $w$ with $\norm{w}_\infty < d$ and $r\in [0,d)$, where $d = \sup_{r \in \R_+}\alpha_1(r)$, define $\hat{\alpha}_4(r) = \max\{ \mathcal{V}(f(\xi,\mu)) | \norm{\mu} \le r, \alpha_1( \mathcal{V}(\xi)) \le r \}$ and $\alpha_4(r) = \max\{\hat{\alpha}_4(r), \alpha_1^{-1}(r)\}$. Then, $\alpha_4$ is continuous, non-decreasing, zero at zero, and can be assumed as a $\mathcal{K}_{[0,d)}$-function (add $\frac{r}{d-r}$ if it is not). Let $\mathcal{S}_w = \{\xi \in \mathcal{S}|\mathcal{V}(\xi) \le \alpha_4(\norm{w}_\infty)\}$. We first show that $\mathcal{S}_w$ is forward invariant.
\\
\textit{Claim:} If $\chi({k_0}) \in \mathcal{S}_w$ for some $k_0 \in \mathbb{Z}_+$, then $\chi(k) \in \mathcal{S}_w$ for all $k \ge k_0$. \\
\textit{Proof of the claim:} Suppose that $\chi({k_0}) \in \mathcal{S}_w$. When $\alpha_1^{-1}(\norm{w}_\infty) < \mathcal{V}(\chi(k_0)) \le \alpha_4(\norm{w}_\infty)$, it follows from \eqref{eq:compaLya} that $\mathcal{V}(\chi(k_0+1)) \le \mathcal{V}(\chi(k_0)) \le \alpha_4(\norm{w}_\infty)$. When $0 \le \mathcal{V}(\chi(k_0)) \le \alpha_1^{-1}(\norm{w}_\infty)$, by the definition of $\hat{\alpha}_4(\norm{w}_\infty)$, $\mathcal{V}(\chi(k_0+1)) \le \hat{\alpha}_4(\norm{w}_\infty) \le {\alpha}_4(\norm{w}_\infty)$. Induction can be used to show that $\mathcal{V}(\chi(k_0+j)) \le \alpha_4(\norm{w}_\infty)$ for all $j \in \mathbb{Z}_+$. 

We now let $k_1 = \min\{ k \in \mathbb{Z}_+| \chi(k) \in \mathcal{S}_w \} \le \infty$. Then, it follows from the above claim that 
\begin{align}\label{eq:ISSSuff1}
    \mathcal{V}(\chi(k)) \le \alpha_4(\norm{w}_\infty), \forall k \ge k_1
\end{align}
For $k < k_1$, $\mathcal{V}(\chi(k)) > \alpha_4(\norm{w}_\infty) \ge \alpha_1^{-1}(\norm{w}_\infty)$, and it follows from \eqref{eq:compaLya} that 
\begin{align}
    \mathcal{V}(\chi(k+1)) - \mathcal{V}(\chi(k)) \le -\alpha_2(\mathcal{V}(\chi(k))).
\end{align}
By the comparison lemma \cite[Lemma 4.3]{Jiang2001ISS}, there exists a $\mathcal{KL}$-function $\beta$ such that 
\begin{align}\label{eq:ISSSuff2}
    \mathcal{V}(\chi(k)) \le \beta(\mathcal{V}(\chi(0)),k), \, \forall k < k_1.
\end{align}
Hence, we can conclude from \eqref{eq:ISSSuff1} and \eqref{eq:ISSSuff2} that system \eqref{eq:nonLinearSys} is small-disturbance ISS with $\gamma = \alpha_4$. 

\textbf{Necessity:} The case when $\mathcal{S} = \Rn$ is proved first. For any bounded inputs $w$ with $\norm{w}_\infty < d$, let $v(k) = \frac{d}{d-\norm{w(k)}}w(k)$. By inverting the function (see Lemma \ref{lm:homomorphism}), we can obtain that $w(k) = \frac{d}{d + \norm{v(k)}}v(k)$ and $\norm{w(k)} = \gamma_1(\norm{v(k)})$, where $\gamma_1(r) = \frac{dr}{d + r} $ is a $\mathcal{K}$-function with the range $[0,d)$. Consider $v(k)$ as the input of system \eqref{eq:nonLinearSys}, we have that
\begin{align}\label{eq:sysRewrite} 
\begin{split}
    \chi(k+1) &= f\left(\chi(k),\frac{d}{d + \norm{v(k)}}v(k)\right) \\
    &= f_1(\chi(k), v(k)).
\end{split}
\end{align}

Since system \eqref{eq:nonLinearSys} is small-disturbance ISS, by \eqref{eq:ISS}, it holds that
\begin{align}\label{eq:sdISSConsequence}
    \mathcal{V}(\chi(k)) \leq \beta(\mathcal{V}(\chi(0)),k) + \gamma \circ \gamma_1(\norm{v}_\infty)
\end{align}
for all $k \in \mathbb{Z}_+$. Since $\mathcal{V}(\xi)$ and $\norm{\xi - \chi^*}$ are size functions for $(\Rn, \chi^*)$, from Lemma \ref{lm:radiallyUnbounded}, there exist $\rho_1, \rho_2 \in \mathcal{K}_\infty$, such that
\begin{align}\label{eq:VsizeNorm}
    \rho_1(\norm{\xi - \chi^*}) \le \mathcal{V}(\xi) \le \rho_2(\norm{\xi - \chi^*}), \, \forall \xi \in \R^n.
\end{align}
The fact that $\gamma \circ \gamma_1 \in \mathcal{K}$, and plugging \eqref{eq:VsizeNorm} into \eqref{eq:sdISSConsequence} and considering Lemma \ref{lm:weakTriangle}, implies that the system \eqref{eq:sysRewrite} is ISS. Consequently, there exists an ISS-Lyapunov function $\mathcal{V}_1$ \cite[Theorem 1]{Jiang2001ISS}, such that 
\begin{align}
    &\rho_3(\norm{\xi- \chi^*} ) \le \mathcal{V}_1(\xi) \le \rho_4(\norm{\xi- \chi^*}), 
\end{align}
and
\begin{align}
    &\mathcal{V}_1(f_1(\xi,\nu)) - \mathcal{V}_1(\xi) \le -\rho_5(\norm{\xi - \chi^*}) + \gamma_2(\norm{\nu}),
\end{align}
for all $ \xi \in \Rn, \nu \in \R^m$, where $\rho_i, \gamma_2 \in \mathcal{K}_\infty \,(i=3,4,5)$. This further implies that, if $\norm{\nu} \le \gamma_2^{-1} \circ \frac{1}{2}\rho_5 \circ \rho_4^{-1}(\mathcal{V}_1(\xi))$, $\mathcal{V}_1(f_1(\xi,\nu)) - \mathcal{V}_1(\xi) \le -\frac{1}{2}\rho_5 \circ \rho_4^{-1}(\mathcal{V}_1(\xi))$. Consequently, if $\mu = \frac{d}{d+\norm{\nu}}\nu$ and $\norm{\mu} \le \gamma_1 \circ \gamma_2^{-1} \circ \frac{1}{2}\rho_5 \circ \rho_4^{-1}(\mathcal{V}_1(\xi))$, $\mathcal{V}_1(f(\xi,\mu)) - \mathcal{V}_1(\xi) \le -\frac{1}{2}\rho_5 \circ \rho_4^{-1}(\mathcal{V}_1(\xi))$. Since $\gamma_1 \circ \gamma_2^{-1} \circ \frac{1}{2}\rho_5 \circ \rho_4^{-1}$ is a $\mathcal{K}$-function with the range $[0,d)$, we conclude that $\mathcal{V}_1$ is a small-disturbance ISS-Lyapunov function for system \eqref{eq:nonLinearSys}.

Next consider the case where $\mathcal{S}$ is an open subset of $\Rn$ which is homeomorphic to $\Rn$. Let $\varphi$ denote the homeomorphism from $\mathcal{S}$ to $\Rn$. For $\chi(k) \in \mathcal{S}$, let $\zeta(k) = \varphi(\chi(k))$ and $\zeta^* = \varphi(\chi^*)$. It follows from \eqref{eq:nonLinearSys} that 
\begin{align}\label{eq:nonlinearSysHomeo}
    \zeta(k+1) = \varphi(f(\varphi^{-1}(\zeta(k)), w(k))) = f_2(\zeta(k), w(k)).
\end{align}
Since system \eqref{eq:nonLinearSys} is small-disturbance ISS over $\mathcal{S} \times \R^m$, it holds that
\begin{align}
    \mathcal{V} \circ \varphi^{-1} (\zeta(k)) \le \beta(\mathcal{V} \circ \varphi^{-1} (\zeta(0)), k ) + \gamma(\norm{w}_\infty) 
\end{align}
for all $k \in \mathbb{Z}_+$ and all $w\in \ell_\infty^{m}$ bounded by $d$. Since the compactness of $\{\xi \in \mathcal{S}| \mathcal{V}(\xi) \le c \}$ implies the compactness of $\{\zeta \in \Rn| \mathcal{V}\circ \varphi^{-1}(\zeta) \le c \}$, $ \mathcal{V}\circ \varphi^{-1}$ is a proper/coercive function over $\Rn$. Hence, system \eqref{eq:nonlinearSysHomeo} is small-disturbance ISS over $\Rn \times \R^m$. According to the results for the case of $\mathcal{S} = \Rn$, there exists a small-disturbance ISS-Lyapunov function $\mathcal{V}_2$ for system \eqref{eq:nonlinearSysHomeo}. It can be easily checked that $\mathcal{V}_2\circ \varphi : \mathcal{S} \to \mathbb{R}_+$ is a small-disturbance ISS-Lyapunov function for system \eqref{eq:nonLinearSys}. In conclusion, the necessity holds.
\end{pf}

In the above analysis, it is required that $\mathcal{S}$ be homeomorphic to $\mathbb{R}^n$, which can often be verified directly from the context of the problem—for example, in the LQR problem, where the set of stabilizing gains is homeomorphic to a Euclidean space. If this condition cannot be verified directly, the arguments in \cite[Theorem 2.2]{WILSON1967323} can be followed to establish that $\mathcal{S}$ is homeomorphic to $\mathbb{R}^n$, provided that the function $f(\cdot, 0)$ is additionally assumed to be a diffeomorphism.

\begin{rem}
    If $f(\cdot ,0)$ is a diffeomorphism and the domain of asymptotic stability of $\chi^*$ is $\mathcal{S}$, then $\mathcal{S}$ is diffeomorphic to $\Rn$. 
\end{rem}

\begin{rem}
\Review{
    We note that small-disturbance ISS is stronger than local ISS. Small-disturbance ISS only requires the disturbance to be sufficiently small, while the initial condition may be arbitrarily far from the equilibrium. By contrast, local ISS requires both the disturbance and the initial condition to be sufficiently small.
}
\Review{
    As an example, consider
    \[
    x_{k+1}=x_k-\frac{2x_k}{1+x_k^2}+w_k .
    \]
    For the nominal system $w_k\equiv 0$, the Lyapunov function
    \[
    V(x)=\ln(1+x^2)
    \]
    satisfies
    \[
    V(x_{k+1})-V(x_k)
    =
    \ln\!\left(1-\frac{4x_k^4}{(1+x_k^2)^3}\right)<0,\quad x_k\neq 0,
    \]
    which shows that the origin is globally asymptotically stable. Hence the system is locally ISS. However, for large $|x_k|$, the stabilizing term $2x_k/(1+x_k^2)$ behaves like $2/x_k$ and becomes arbitrarily small. Therefore, a sufficiently small constant disturbance can dominate this restoring effect when the state is large, implying that the system is not small-disturbance ISS.
}
\end{rem}

\section{Robustness Analysis of Perturbed Gradient Descent}
This section applies the concept of small-disturbance ISS to analyze the gradient descent algorithm for solving the constrained nonlinear program:
\begin{align}\label{eq:optiPro}
    &\min_{z \in \mathcal{Z}} \mathcal{J}(z) 
\end{align}    
where $\mathcal{Z}$ is an admissible set, defined as an open subset of $\mathbb{R}^n$ that is homeomorphic to $\mathbb{R}^n$, and $\mathcal{J}: \mathcal{Z} \to \mathbb{R}$ is an objective function with a unique global minimizer $z^*$. 

\begin{defn}\label{def:ProperLossFun}
    A continuously differentiable function $\mathcal{J}: \mathcal{Z} \to \mathbb{R}$ is a proper objective function if
    \begin{enumerate}
        \item $\mathcal{J}(z) - \mathcal{J}(z^*)$ is a size function for $(\mathcal{Z}, z^*)$;
        \item $\nabla \mathcal{J}(z)$ is Lipschitz continuous over the sublevel set $\mathcal{Z}(h) = \{z \in \mathcal{Z}| \mathcal{J}(z) \le h \}$, that is 
        \begin{align}
            \norm{\nabla \mathcal{J}(z_1) - \nabla\mathcal{J}(z_2)} \le L(h)\norm{z_1 - z_2}
        \end{align}
        for all $z_1,z_2 \in \mathcal{Z}(h)$.
        \item there exists a $\mathcal{K}$-function $\alpha_5$, such that $\norm{\nabla \mathcal{J}(z)} \ge \alpha_5(\mathcal{J}(z) - \mathcal{J}(z^*))$ ($\mathcal{K}$-PL estimate).
    \end{enumerate}
\end{defn}

The perturbed gradient descent method for \eqref{eq:optiPro} is 
\begin{align}\label{eq:gradientDescent}
    {z}(k+1) = z(k) - \eta(k) ( \nabla \mathcal{J}(z(k)) + e(k)),
\end{align}
where $\eta(k) > 0$ is a step size, and $e \in \ell_\infty^{n}$ denotes the perturbation to the gradient algorithm. The perturbation $ e(k)$ may arise from inaccurate gradient estimation in data-driven optimization, rounding errors in numerical computation, or even malicious attacks on the gradient descent algorithm. 

\begin{thm}\label{thm:gradDescentISS}
    If $\mathcal{J}$ is a proper objective function and the step size satisfies $0 < \eta_{\min} \le \eta(k) \le \frac{1}{L(\mathcal{J}(z(k)))}$, then system \eqref{eq:gradientDescent} is small-disturbance ISS.
\end{thm}
\begin{pf}
Denote 
\begin{align}
    \kappa(k,s) = \mathcal{J}\left(z(k) - s( \nabla \mathcal{J}(z(k)) + e(k)) \right). 
\end{align}
The derivative of $\kappa(k,s)$ is written as
\begin{align}\label{eq:kappaDerivative}
\begin{split}
    &\frac{\partial \kappa(k,s)}{\partial s}\bigg\lvert_{s=0} = -\innprod{\nabla \mathcal{J}\left(z(k)\right)}{\nabla \mathcal{J}(z(k)) + e(k)} \\
    &\quad \le -\frac{1}{2}\norm{\nabla \mathcal{J}(z(k))}^2 + \frac{1}{2}\norm{e(k)}^2 \\
    &\quad \le - \frac{1}{2} \alpha_5 \left(\mathcal{J}(z(k)) - \mathcal{J}(z^*) \right)^2  + \frac{1}{2}\norm{e(k)}^2
\end{split}
\end{align}
where the second line is from Cauchy-Schwarz inequality and Young’s inequality, and the last line is a direct consequence of the $\mathcal{K}$-PL property. 

When $\norm{e(k)} \le \frac{1}{2}\alpha_5(\mathcal{J}(z(k)) - \mathcal{J}(z^*))$, we first show that $z(k+1) \in \mathcal{Z}(\mathcal{J}(z(k)))$. It follows from \eqref{eq:kappaDerivative} that 
\begin{align}
    \frac{\partial \kappa(k,s)}{\partial s}\bigg\lvert_{s=0} \le - \frac{3}{8} \alpha_5 \left(\mathcal{J}(z(k)) - \mathcal{J}(z^*) \right)^2 < 0.
\end{align}
Hence, $\kappa(k,s) < \kappa(k,0) = \mathcal{J}(z(k))$ for some small $s>0$. Since $\mathcal{J}$ is coercive, the sublevel set $\mathcal{Z}(\mathcal{J}(z(k)))$ is compact \cite[Lemma 2.4]{Sontag2022}. As a result, $z(k) - s( \nabla \mathcal{J}(z(k)) + e(k))$ can reach the boundary of $\mathcal{Z}(\mathcal{J}(z(k)))$ for some $\bar{s}(k)>0$. Denote $\bar{s}(k)$ as the point at which the boundary is first reached, i.e. $\kappa(k,\bar{s}(k)) = \mathcal{J}(z(k))$, and $z(k) - s( \nabla \mathcal{J}(z(k)) + e(k)) \in \mathcal{Z}(\mathcal{J}(z(k)))$ for all $0 \le s\le \bar{s}(k)$. Since $\kappa(k,s)$ is $L_1(\mathcal{J}(z(k)))$-smooth over $s \in [0,\bar{s}(k)]$ with $L_1(\mathcal{J}(z(k))) = \norm{\nabla \mathcal{J}(z(k))+e(k)}^2 L(\mathcal{J}(z(k)))$, according to \cite[Lemma 1.2.3]{nesterov2013introductory}, it holds that
\begin{align}\label{eq:kappaDerivative1}
\begin{split}
   &\kappa(k,\bar{s}(k)) \le \kappa(k,0) + \frac{\partial \kappa(k,s)}{\partial s}\bigg\lvert_{s=0}\bar{s}(k)  \\
   &\qquad + \frac{L(\mathcal{J}(z(k)))}{2}\norm{\nabla \mathcal{J}(z(k))+e(k)}^2\bar{s}(k)^2 \\
   &\quad = \kappa(k,0) - (\bar{s}(k) - \frac{L(\mathcal{J}(z(k)))}{2}\bar{s}(k)^2) \norm{\nabla \mathcal{J}(z(k))}^2 \\
   &\qquad - (\bar{s}(k) - L(\mathcal{J}(z(k)))\bar{s}(k)^2) \innprod{\nabla \mathcal{J}(z(k))}{e(k)} \\
   &\qquad + \frac{L(\mathcal{J}(z(k)))}{2}\bar{s}(k)^2\norm{e(k)}^2
\end{split}
\end{align}
where the second line is from \eqref{eq:kappaDerivative}. Suppose that $\bar{s}(k) < \frac{1}{L(\mathcal{J}(z(k)))}$, which implies
\begin{align}
(\bar{s}(k) - L(\mathcal{J}(z(k)))\bar{s}(k)^2) > 0.
\end{align}
Using Cauchy-Schwarz inequality, Young's inequality, and the $\mathcal{K}$-PL condition, 
\begin{align}\label{eq:kappaDerivative2}
\begin{split}
&\kappa(k,\bar{s}(k)) \le \kappa(k,0) \\
&\quad - \frac{\bar{s}(k)}{2} \alpha_5 \left(\mathcal{J}(z(k)) - \mathcal{J}(z^*) \right)^2  + \frac{\bar{s}(k)}{2} \norm{e(k)}^2\\
&\le \kappa(k,0) - \frac{3\bar{s}(k)}{8} \alpha_5 \left(\mathcal{J}(z(k)) - \mathcal{J}(z^*) \right)^2   \\
&< \mathcal{J}(z(k)). 
\end{split}
\end{align}
Since $\kappa(k,\bar{s}(k)) < \mathcal{J}(z(k))$ contradicts $\kappa(k,\bar{s}(k)) = \mathcal{J}(z(k))$, it follows that $\bar{s}(k) \ge \frac{1}{L(\mathcal{J}(z(k)))}$. Note that $z(k) - s( \nabla \mathcal{J}(z(k)) + e(k))$ reaches the boundary of $\mathcal{Z}(\mathcal{J}(z(k)))$ for the first time when $s = \bar{s}(k)$. Hence, $\eta(k) \le\frac{1}{L(\mathcal{J}(z(k)))} \le  \bar{s}(k)$ leads to $z(k+1) \in \mathcal{Z}(\mathcal{J}(z(k)))$.

Next, we show that $\mathcal{J}(z(k)) - \mathcal{J}(z^*)$ is a small-disturbance ISS-Lyapunov function for system \eqref{eq:gradientDescent}. Note that
\begin{align}
    \mathcal{J}(z(k+1)) = \kappa(k, \eta(k)).
\end{align}
Since $z(k), z(k+1) \in \mathcal{Z}(\mathcal{J}(z(k)))$ when $\norm{e(k)} \le \frac{1}{2}\alpha_5(\mathcal{J}(z(k)) - \mathcal{J}(z^*))$, by \eqref{eq:kappaDerivative1} and \eqref{eq:kappaDerivative2}, it holds that
\begin{subequations}
\begin{align}
     &\mathcal{J}(z(k+1)) - \mathcal{J}(z(k)) \le - \frac{\eta(k)}{2} \alpha_5 \left(\mathcal{J}(z(k)) - \mathcal{J}(z^*) \right)^2 \nonumber\\
     &\qquad + \frac{\eta(k)}{2} \norm{e(k)}^2  \label{eq:Jdiff}\\
     &\quad \le - \frac{3\eta(k)}{8} \alpha_5 \left( \mathcal{J}(z(k)) - \mathcal{J}(z^*) \right)^2 . \label{eq:Jdecrease}
\end{align}
\end{subequations}
According to Theorem \ref{thm:smallISSSufficient}, we conclude that system \eqref{eq:gradientDescent} is small-disturbance ISS.
\end{pf}

In conclusion, Theorem \ref{thm:gradDescentISS} indicates that despite the presence of perturbation $e$ with a bounded magnitude $\norm{e}_\infty < \frac{1}{2}\sup_{r \in \R_+}\alpha_5(r)$, the gradient descent algorithm can still converge to a neighborhood of the optimum $z^*$, specifically $\{z \in \mathcal{Z}| \mathcal{J}(z) - \mathcal{J}(z^*) \le \alpha_5^{-1}(2\norm{e}_\infty)\}$. 

Corollary \ref{cor:ISSgradient} indicates that if the $\mathcal{K}$-PL estimate in Definition \ref{def:ProperLossFun} is strengthened to a $\mathcal{K}_\infty$-function, then the system \eqref{eq:gradientDescent} is ISS \cite{Jiang2001ISS}, which is a stronger property than small-disturbance ISS.

\begin{cor}\label{cor:ISSgradient}
    If the $\mathcal{K}$-PL estimate $\alpha_5$ in Definition \ref{def:ProperLossFun} is strengthened to a $\mathcal{K}_\infty$-function, and the step size satisfies $0 < \eta_{\min} \le \eta(k) \le \frac{1}{L(\mathcal{J}(z(k)))}$, then system \eqref{eq:gradientDescent} is ISS.    
\end{cor}
\begin{pf}
    By following the same reasoning as in the proof of Theorem \ref{thm:gradDescentISS}, if $\norm{e(k)} \le \frac{1}{2}\alpha_5(\mathcal{J}(z(k)) - \mathcal{J}(z^*))$, then
    \begin{align}
    \begin{split}
         &\mathcal{J}(z(k+1)) - \mathcal{J}(z(k)) \\
         &\quad \le - \frac{3\eta(k)}{8} \alpha_5 \left( \mathcal{J}(z(k)) - \mathcal{J}(z^*) \right)^2. 
    \end{split}
    \end{align}
    Since $\alpha_5$ is a $\mathcal{K}_\infty$-function, it follows from \cite[Lemma 3.5]{Jiang2001ISS} that the system \eqref{eq:gradientDescent} is ISS.    
\end{pf}

Additionally, the following corollary states that, if the $\mathcal{K}$-PL estimate $\alpha_5$ in Definition \ref{def:ProperLossFun} is relaxed to a positive definite function, then the system \eqref{eq:gradientDescent} is integral ISS, which is equivalent to global asymptotic stability for discrete-time dynamical systems \cite{angeli1999intrinsic}.

\begin{cor}\label{cor:GradientiISS}
    If the $\mathcal{K}$-PL estimate $\alpha_5$ in Definition \ref{def:ProperLossFun} is relaxed to a positive definite function, and the step size satisfies $0 < \eta_{\min} \le \eta(k) \le \frac{1}{L(\mathcal{J}(z(k)))}$, then system \eqref{eq:gradientDescent} is integral ISS.     
\end{cor}
\begin{pf}
    In the absence of noise, the gradient system \eqref{eq:gradientDescent} satisfies the inequality:
    \begin{align}
    \begin{split}
         &\mathcal{J}(z(k+1)) - \mathcal{J}(z(k)) \\
         &\quad \le - \frac{3\eta(k)}{8} \alpha_5 \left( \mathcal{J}(z(k)) - \mathcal{J}(z^*) \right)^2. 
    \end{split}
    \end{align}    
    This implies that the gradient system is globally asymptotically stable. Moreover, by the main result of \cite{angeli1999intrinsic}, it is also integral ISS.
\end{pf}

\begin{rem}\label{rem:ISSsdISS}
    The gradient dominance condition proposed in Definition \ref{def:ProperLossFun}, i.e. the $\mathcal{K}$-PL condition, can be viewed as a nonlinear generalization of the well-known PL condition. If the classical PL condition holds \cite{Polyak1963,Lojasiewicz1963}, meaning $\alpha_5(r) = c \sqrt{r}$ for all $r \ge 0$ and some $c>0$, the perturbed gradient descent algorithm in \eqref{eq:gradientDescent} is exponentially ISS.
\end{rem}

\begin{rem}
\Review{
    If the perturbations are state-dependent and satisfy
    \[
    \|e(z(k))\|^2 \le \frac{1}{4}\alpha_5\big(\mathcal{J}(z(k))-\mathcal{J}(z^*)\big)^2 + r(k), \quad \forall k \in \mathbb{Z}_+,
    \]
    then, by \eqref{eq:Jdiff}, the state-dependent term can be absorbed into the Lyapunov analysis, while $r(k)$ acts as an external disturbance. Consequently, the gradient dynamics remains small-disturbance ISS with respect to $r(k)$.}
\end{rem}

\Review{
We introduce a structured $\mathcal{K}$-PL condition, referred to as the $\mathcal{K}_{\mathrm{SAT}}$-PL condition \cite{Arthur2025}, which specifies the comparison function $\alpha_5$ as
\begin{align}\label{eq:KSAT}
\alpha_5(r) = \sqrt{\frac{\theta_1 r}{r + \theta_2}}, 
\qquad \forall r \ge 0,
\end{align}
for some constants $\theta_1>0$ and $\theta_2>0$. This form captures a saturation effect for large $r$ while preserving the classic PL condition near the origin. Under the $\mathcal{K}_{\mathrm{SAT}}$-PL condition, for any $\bar r>0$, one may choose $\theta_{3,\bar r} = \sqrt{{\theta_1}/{(\bar r + \theta_2)}}$
such that
\begin{align}\label{eq:KSATPLLinear}
\alpha_5(r) \ge \theta_{3,\bar r}\sqrt r,
\qquad \forall r \in [0,\bar r].
\end{align}
Let $\bar r \ge 
\max\!\left\{
\mathcal J(z(0))-\mathcal J^*,
\,
\alpha_5^{-1}(2\|e\|_\infty)
\right\}$. Whenever
\[
\|e(k)\|
\le
\tfrac12
\alpha_5\!\big(
\mathcal J(z(k))-\mathcal J(z^*)
\big),
\]
inequalities \eqref{eq:Jdecrease} and \eqref{eq:KSATPLLinear} imply
\begin{align}
\mathcal J(z(k+1))-\mathcal J(z(k))
\le
-\frac{3\eta(k)}{8}\theta_{3,\bar r}^2
\big(
\mathcal J(z(k))-\mathcal J^*
\big).
\end{align}
Therefore, whenever
\[
\mathcal J(z(k))-\mathcal J(z^*)
\ge
\alpha_5^{-1}(2\|e\|_\infty),
\]
the optimality gap satisfies the linear convergence
\begin{align*}
&\mathcal J(z(k+1))-\mathcal J(z^*) \\
& \quad \le\left(
1-
\frac{3\eta(k)}{8}\theta_{3,\bar r}^2
\right)
\big(
\mathcal J(z(k))-\mathcal J(z^*)
\big).
\end{align*}
Since $\lim_{\bar r \to \infty} \theta_{3,\bar r} = 0$, the guaranteed linear convergence rate deteriorates as the initial condition moves farther from the minimizer. Consequently, the $\mathcal{K}_{\mathrm{SAT}}$-PL condition implies not only the small-disturbance ISS property of the perturbed gradient descent dynamics, but also a semi-global linear convergence rate within any bounded sublevel set.
}

We present several examples of objective functions for which the robustness of the associated gradient descent algorithms is analyzed using the proposed framework.

\begin{exmp}
Consider the optimization problem $\min_{z \in (1,\infty)} \mathcal{J}(z) = \frac{(z-z^*)^2}{2(z-1)}$ where $z^* = 1+\sqrt{2}$ and $\mathcal{J}(z^*) = 0$. It can be verified that
\begin{align*}
    \| \nabla\mathcal{J}(z) \| \ge \frac{\mathcal{J}(z)}{\sqrt{2}/2 + 2\mathcal{J}(z)}.
\end{align*}
The function $\mathcal{J}(z)$ is a size function and $\nabla\mathcal{J}(z)$ is Lipschitz continuous over any sublevel sets. Consequently, if the step size satisfies $0 < \eta(k) \le \frac{1}{L(\mathcal{J}(z(k)))}$, then the perturbed gradient descent algorithm in \eqref{eq:gradientDescent} is small-disturbance ISS. 
\end{exmp}

\begin{exmp}
Consider the optimization problem $\min_{z \in \R}\mathcal{J}(z) = \frac{1}{4}z^4$ with the optimal solution $z^* = 0$ and $\mathcal{J}(z^*) = 0$. The function is not strongly convex since its second derivative $\nabla^2 \mathcal{J}(z) = 0$ at $z = 0$. Furthermore, since $\norm{\nabla \mathcal{J}(z)}^2 = o(\mathcal{J}(z))$ as $z \to 0$, there does not exist a constant $c>0$ such that $\norm{\nabla \mathcal{J}(z)}^2 \ge c \mathcal{J}(z)$, which implies that the classical PL condition does not hold. Since $\norm{\nabla \mathcal{J}(z)} = (4\mathcal{J}(z))^{{3}/{4}}$, it follows from Corollary \ref{cor:ISSgradient} that the perturbed gradient descent in \eqref{eq:gradientDescent} is ISS if the step size satisfies $0 < \eta(k) \le \frac{1}{L(\mathcal{J}(z(k)))}$. 
\end{exmp}

\begin{exmp}
    Consider the optimization problem $\min_{z \in \R} \mathcal{J}(z) = \log(z^2+1)$, which has a minimum at $z^* = 0$. The gradient is given by $\nabla \mathcal{J}(z) = \frac{2z}{z^2+1}$. Since $\lim_{z \to \infty} \norm{\nabla \mathcal{J}(z)} = 0$ while $\lim_{z \to \infty} \mathcal{J}(z) = \infty$, we can only identify a positive definite function $\alpha_5$ such that $\norm{\nabla \mathcal{J}(z)} \ge \alpha_5(\mathcal{J}(z))$. Therefore, by Corollary \ref{cor:GradientiISS}, the gradient system \eqref{eq:gradientDescent} is integral ISS provided that the step size satisfies $0 < \eta(k) \le \frac{1}{L(\mathcal{J}(z(k)))}$. 
\end{exmp}


\section{Application to LQR Problem}
In this section, we utilize the tool of small-disturbance ISS to analyze the robustness of the gradient descent algorithms in solving the LQR problem. Some preliminaries on the LQR are introduced in the next subsection.
\subsection{Preliminaries of LQR}
A linear time-invariant (LTI) system can be represented by 
\begin{align}\label{eq:LTI}
    \dot{x}(t) = Ax(t) + Bu(t), \,\, x(0) = x_0
\end{align}
where $x(t) \in \mathbb{R}^n$ denotes the state with the initial state $x_0$; $u(t) \in \mathbb{R}^{m}$ denotes the control input; $A \in \R^{n\times n}$ and $B \in \mathbb{R}^{n\times m}$ are constant matrices. The LQR problem aims to find a state-feedback controller by minimizing the cumulative quadratic cost
\begin{align}\label{eq:LQRcost}
    \mathcal{J}_{1}(x_0, u) =  \int_{0}^{\infty} x(t)^\top Qx(t) + u(t)^\top R u(t) \de t , 
\end{align}
with $Q \in \mathbb{S}^n_{++}$ and $R \in \mathbb{S}^m_{++}$. Suppose that $(A,B)$ is stabilizable, according to \cite[Section 8.4]{book_sontag}, the optimal controller is expressed as
\begin{align}\label{eq:Koptexpression}
    u^*(x(t)) = -{K^*}x(t), \quad K^* = R^{-1}B^\top P^*,
\end{align}
where $P^* \in \mathbb{S}_{++}^n$ is the solution of the algebraic Riccati equation (ARE)
\begin{align}\label{eq:continuousARE}
    A^\top  P^* + P^* A + Q - P^*BR^{-1}B^\top P^* = 0.
\end{align}

The optimal controller $K^*$, as indicated in \cite{CJS2024}, can be found by directly solving the optimization problem 
\begin{align}\label{eq:costJc_closedform}
    \min_{K \in \mathcal{G}}\mathcal{J}_2(K) := \Tr{P_K},
\end{align}
where $\mathcal{G} = \{K \in \mathbb{R}^{m \times n}| A-BK \text{ is Hurwitz} \}$ is the admissible set containing all the stabilizing gains, and $P_K \in \mathbb{S}^n_{++}$ is the solution of the Lyapunov equation
\begin{align}\label{eq:continuousLyapunov}
    (A-BK)^{\top}P_K + P_K (A-BK) + Q + K^\top  R K = 0.
\end{align}
Since $\mathcal{G}$ is diffeomorphic to an open convex subset of $\R^{m \times n}$ \cite[Lemma 3.2]{bu2019topological}, and every open convex subset of $\R^{m \times n}$ is homeomorphic to $\R^{m \times n}$ itself \cite{geschke2012convex}, it follows that $\mathcal{G}$ is  homeomorphic to $\R^{m \times n}$.

We next derive the gradient of $\mathcal{J}_2(K)$. For any $K \in \mathcal{G}$, the increment of \eqref{eq:continuousLyapunov} is
\begin{align}\label{eq:PKderivative}
\begin{split}
    &(A-BK)^{\top}\de P_K + \de P_K (A-BK) \\
    & \quad + \de K^\top(RK - B^\top P_K)+ (RK - B^\top P_K)^\top \de K  = 0.    
\end{split}
\end{align}
It follows from Corollary \ref{cor:traceLya} that
\begin{align}\label{eq:dPKdK}
    \Tr{\de P_K} = 2\innprod{\de K}{(RK - B^\top P_K) Y_K}
\end{align}
where $Y_K \in \mathbb{S}^n_{++}$ is the solution of 
\begin{align}\label{eq:YKDef}
    (A-BK)Y_K + Y_K (A-BK)^\top  + I_n = 0.
\end{align}
Hence, for any $K \in \mathcal{G}$, the gradient of $\mathcal{J}_2(K)$ can be computed by 
\begin{align}\label{eq:gradientLQR}
    \nabla \mathcal{J}_2(K) = 2(RK - B^\top P_K)Y_K.
\end{align}

Throughout this article, define $Y^* \in \mathbb{S}_{++}^n$ as the solution of \eqref{eq:YKDef} with $K$ replaced by $K^*$. In addition, using corollary \ref{cor:traceLya} and Lemma \ref{lm:traceIneq}, it holds that
\begin{align}\label{eq:traceYbound}
    \mathcal{J}_2(K) = \Tr{(Q+K^\top R K)Y_K} \ge \eigmin{Q} \Tr{Y_K}.
\end{align}

The following lemma shows that the gradient $\nabla \mathcal{J}_2(K)$ is Lipschitz continuous over the sublevel set $\mathcal{G}(h) = \{K \in \mathcal{G}| \mathcal{J}_2(K) \le h\}$. 
\begin{lem}[Theorem 3.15 in \cite{Fatkhullin2021}] \label{lm:Lsmoothness}
    The gradient $\nabla \mathcal{J}_2(K)$ is $L(h)$-Lipschitz continuous over the sublevel set $\mathcal{G}(h) = \{K \in \mathcal{G}| \mathcal{J}_2(K) \le h\}$.
\end{lem}

\subsection{Small-Disturbance ISS of Standard Gradient Descent}
This subsection applies the concept of small-disturbance ISS to analyze the robustness of the standard gradient descent method for the LQR problem,
\begin{align}\label{eq:graddescentLQR}
\begin{split}
    &K(k+1) = K(k) -\eta(k) (\nabla \mathcal{J}_2(K(k)) + W(k)) \\
    & = K(k)  -\eta(k) (2(RK(k) - B^\top P(k))Y(k) + W(k)), 
\end{split}
\end{align}
where $P(k) = P_{K(k)}$, $Y(k) = Y_{K(k)}$, and $W \in \ell_{\infty}^{m \times n}$ is the perturbation to the gradient descent algorithm. The perturbation $W$ can represent gradient estimation errors in the context of data-driven control. When the system matrices are unknown, gradient estimation can be achieved through the finite-difference method \cite{fazel2018global} or approximate dynamic programming \cite{tutorial_Jiang}, both of which introduce errors due to measurement noise, system process noise, and even potential malicious attacks on the algorithm. The following lemma is introduced to ensure that $\mathcal{J}_2(K)$ satisfyies the $\mathcal{K}$-PL condition in Definition \ref{def:ProperLossFun}, which is critical to the robustness analysis.

\begin{lem} [Lemma 5.7 in \cite{CJS2024}] \label{lm:gradientClassK}
The objective function $\mathcal{J}_2(K)$ satisfies the $\mathcal{K}$-PL condition, that is 
\begin{align*}
    \norm{\nabla \mathcal{J}_2(K)}_F \geq  \alpha_6(\mathcal{J}_2(K) - \mathcal{J}_2(K^*)), \quad \forall K \in \mathcal{G},
\end{align*}
where 
\begin{align}
     \alpha_6(r) = \frac{r}{b_1r + b_2},
\end{align}
\begin{align*}
b_1 = \frac{\norm{B}\sqrt{2(\eigmin{Y^*}+\eigmax{Y^*})}}{\eigmin{R} \sqrt{\eigmin{Y^*}}}, 
\end{align*}
and 
\begin{align*}
&b_2 = \\
&\frac{\norm{A-BK^*}_F^2\eigmin{Y^*}^{{1}/{2}}(\eigmin{Y^*}+\eigmax{Y^*})^{{1}/{2}}}{\sqrt{2}\norm{B}}
\end{align*}
\end{lem}
\Review{By Lemma \ref{lm:gradientClassK}, \cite[Corollary 3.5.1]{bu2020policy}, and \cite[Lemma 5]{Arthur2025}, the LQR cost $\mathcal{J}_2$ satisfies the $\mathcal{K}_{\mathrm{SAT}}$-PL condition. Therefore, the gradient descent iteration associated with $\mathcal{J}_2$ achieves semi-global linear convergence.}

\begin{rem}
The classical PL condition requires $\norm{\nabla \mathcal{J}_2(K)}^2 \ge c (\mathcal{J}_2(K) - \mathcal{J}_2(K^*))$, which holds only over a compact sublevel set $\mathcal{G}(h) = \{K \in \mathcal{G}| \mathcal{J}_2(K) \le h\}$ due to the limitations of its linear form \cite{bu2020policy,Mohammadi2022}. In this article, by generalizing the classical PL condition to nonlinear form, we can obtain a global estimate of the gradient dominance condition.
\end{rem}

With the established Lipschitz continuity and the $\mathcal{K}$-PL condition for the objective function $\mathcal{J}_2(K)$, we are now prepared to present the main result on the robustness of the gradient descent algorithm in \eqref{eq:graddescentLQR}.
\begin{thm}\label{thm:LQRgradientISS}
    If $0 < \eta_{\min} \le \eta(k) \le \frac{1}{L(\mathcal{J}_2(K(k)))}$, where $L(h)$ is defined in Lemma \ref{lm:Lsmoothness}, the standard gradient descent algorithm in \eqref{eq:graddescentLQR} is small-disturbance ISS with respect to $W$.  
\end{thm}
\begin{pf}
Our task is simply to prove that the properties in Definition \ref{def:ProperLossFun} are satisfied by the objective function $\mathcal{J}_2(K)$. The coercivity of $\mathcal{J}_2(K)$ is demonstrated in \cite[Lemma 3.3]{bu2020policy}. The $L(h)$-Lipschitz continuity and the $\mathcal{K}$-PL condition are established in Lemmas \ref{lm:Lsmoothness} and \ref{lm:gradientClassK}, respectively. Consequently, by Theorem \ref{thm:gradDescentISS}, the gradient descent algorithm in \eqref{eq:graddescentLQR} is small-disturbance ISS. 
\end{pf}
\subsection{Small-Disturbance ISS of Natural Gradient Descent}
This subsection analyzes the robustness of the natural gradient descent algorithm, developed by leveraging the Riemannian geometry of the objective function $\mathcal{J}_2(K)$. By subtracting \eqref{eq:continuousARE} from \eqref{eq:continuousLyapunov} and completing the squares, the Lyapunov equation
\begin{align}\label{eq:PPoptLyap}
\begin{split}
    &(A-BK)^\top (P_K - P^*) + (P_K - P^*) (A-BK) \\
    &\quad + (K-K^*)^\top R (K-K^*) = 0   
\end{split}
\end{align}
can be obtained. Applying Corollary \ref{cor:traceLya} to \eqref{eq:PPoptLyap}, the LQR cost can be expressed as a quadratic function over the Riemannian metric $(\mathcal{G}, \innprod{\cdot}{\cdot}_{Y_K})$, that is, 
\begin{align*}
    \mathcal{J}_2(K) - \mathcal{J}_2(K^*) = \innprod{K-K^*}{R(K-K^*)}_{Y_K}.
\end{align*}

The standard gradient descent in \eqref{eq:graddescentLQR} follows the steepest descent direction under the standard Euclidean metric $(\mathcal{G}, \innprod{\cdot}{\cdot}_{I_n})$.  However, this ad hoc choice of metric may not be appropriate.  As seen in the expression for $\nabla \mathcal{J}_2(K)$ in \eqref{eq:gradientLQR}, the magnitude of the gradient depends on $Y_K$, which can diverge as $K \to \partial \mathcal{G}$ but vanishes as $\norm{K}_F \to \infty$ (see Example \ref{example:1dimsystem}). The non-isotropic property induced by the improper choice of the Euclidean metric may degrade the convergence rate. As pointed out by Amari \cite{Amari1998,Amari1998Why}, the choice of a metric should be based on the manifold that the optimization parameters lie in. Over the Riemannian manifold $(\mathcal{G}, \innprod{\cdot}{\cdot}_{Y_K})$ and according to \cite{Amari1998Why,CJS2024}, the steepest-descent direction can be derived as
\begin{align*}
    \grad{\mathcal{J}_2(K)} = \nabla \mathcal{J}_2(K) Y_K^{-1} = 2(RK - B^\top P_K).
\end{align*}
In practice, the accurate gradient is not accessible and should be estimated through sampling and experiments. The perturbed natural gradient descent algorithm is 
\begin{align}\label{eq:Gradient_NaturalCont}
    K(k+1) = K(k) - \eta(k) (2(RK(k) - B^\top P(k)) + W(k))
\end{align}
where $P(k) = P_{K(k)}$, $W \in \ell_{\infty}^{m \times n}$ denotes the perturbation and $\eta(k)>0$ is the step size to be determined later.

The following example illustrates the advantage of natural gradient descent over standard gradient descent in terms of convergence rate.
\begin{exmp}\label{example:1dimsystem}
    Consider an LTI system with $A=0$ and $B=Q=R=1$, the admissible set of stabilizing gains is $\mathcal{G} = (0,\infty)$, with $Y_K = \frac{1}{2K}$ and cost function $\mathcal{J}_2(K) = P_K = \frac{K^2 + 1}{2K}$. The standard gradient is $\nabla \mathcal{J}_2(K) = \frac{K^2-1}{2K^2}$, while the natural gradient is $\grad{\mathcal{J}_2(K)} = \frac{K^2-1}{K}$. The optimal solution is $K^* = 1$ with the corresponding optimal cost $\mathcal{J}_2(K^*) = 1$. 

    As $K \to \infty$, the standard gradient $\nabla \mathcal{J}_2(K)$ saturates, while the natural gradient $\grad{\mathcal{J}_2(K)}$ remains unbounded. This distinction allows the natural gradient to achieve faster convergence, particularly when $K$ is far from the optimum. Under standard gradient descent with the update rule
    $$
    K' = K - \eta\frac{K^2-1}{2K^2}, 
    $$
    it can be shown that 
    $$
    \mathcal{J}_2(K') - \mathcal{J}_2(K) = - \eta m_1(K,\eta)  (\mathcal{J}_2(K) - \mathcal{J}_2(K^*)),
    $$
    where $m_1(K,\eta) = \frac{(2K-\eta)(K+1)^2}{4K^4 - 2\eta K^3 + 2\eta K}$.
    In comparison, under natural gradient descent with the update rule
    $$
    K' = K - \eta\frac{K^2 - 1}{K},
    $$
    it can be verified that
    $$
    \mathcal{J}_2(K') - \mathcal{J}_2(K) = - \eta m_2(K,\eta)  (\mathcal{J}_2(K) - \mathcal{J}_2(K^*)),    
    $$
    where $m_2(K,\eta) = \frac{(1-\eta)(K+1)^2}{(1-\eta)K^2 + \eta}$. Since $\lim_{K \to \infty}m_1(K, \eta) = 0$  while $\lim_{K \to \infty}m_2(K, \eta) = 1$, the convergence rate of natural gradient descent is faster than that of standard gradient descent when $K$ is far from the optimum.

\end{exmp}

The following two lemmas are introduced to assist in the development of the small-disturbance ISS property of natural gradient descent. 

\begin{lem}\label{lm:PdiffNatural}
    For any $K\in\mathcal{G}$, the relation \begin{align}\label{eq:KKoptK'Geo}
    \begin{split}
        &2\innprod{K-K^*}{R(K-K_+)}_{Y^*} = \Tr{P_K - P^*} \\
        &\quad + \innprod{K-K^*}{R(K-K^*)}_{Y^*}
    \end{split}
    \end{align}
    holds. Note that $K_+ = R^{-1}B^\top  P_K$.
\end{lem}
\begin{pf}
    Subtracting \eqref{eq:continuousARE} from \eqref{eq:continuousLyapunov} and completing the squares yield    \begin{align}\label{eq:AoptPKPoptDiff}
    \begin{split}
        &(A-BK^*)^\top (P_K - P^*) + (P_K - P^*)(A-BK^*) \\
        &\quad + (K-K_+)^\top R (K-K_+) \\
        &\quad - (K_+ - K^*)^\top R (K_+ - K^*) = 0.
    \end{split}
    \end{align}
    Completing the squares again gives that
    \begin{align}\label{eq:KK'CompletingSquare}
     \begin{split}
         &(K - K_+)^\top  R (K - K_+)  - (K_+ - K^*)^\top  R (K_+ - K^*)  \\
         &\quad = (K - K_+)^\top  R (K - K^*) + (K - K^*)^\top  R (K - K_+) \\
         &\qquad - (K - K^*)^\top  R (K - K^*).  
    \end{split}   
    \end{align}    
    Plugging \eqref{eq:KK'CompletingSquare} into \eqref{eq:AoptPKPoptDiff} results in 
    \begin{align*}
    \begin{split}
        &(A-BK^*)^\top (P_K - P^*) + (P_K - P^*)(A-BK^*) \\
        &\quad + (K - K_+)^\top  R (K - K^*)  +(K - K^*)^\top  R (K - K_+) \\
        &\quad - (K - K^*)^\top  R (K - K^*) = 0. 
    \end{split}
    \end{align*}
    Since $A-BK^*$ is Hurwitz, according to Corollary \ref{cor:traceLya}, it holds that
    \begin{align}
    \begin{split}
        \Tr{P_K - P^*} &= 2\innprod{K-K^*}{R(K-K_+)}_{Y^*} \\
        &\quad - \innprod{K-K^*}{R(K-K^*)}_{Y^*}        
    \end{split}
    \end{align}
    which completes the proof.
\end{pf}

\begin{lem}\label{lm:KK'bound}
    For any $K \in \mathcal{G}$, it holds that
    \begin{align}\label{eq:KK'YoptUppderBound}
        \innprod{K-K_+}{R(K-K_+)}_{Y^*} \le c(K) \Tr{P_K - P^*}
    \end{align}
    where $c(K) = 1+\norm{Y^*}\norm{BR^{-1}B^\top  }\Tr{P_K}$.
\end{lem}
\begin{pf}
    It follows from \eqref{eq:AoptPKPoptDiff} and Corollary \ref{cor:traceLya} that 
    \begin{align}
    \begin{split}
        &\innprod{K-K_+}{R(K-K_+)}_{Y^*} = \Tr{P_K - P^*} \\
        &\quad + \innprod{K_+ - K^*}{R(K_+ - K^*)}_{Y^*} \\
        &\le \Tr{P_K - P^*} + \norm{Y^*}\norm{BR^{-1}B^\top  }\Tr{P_K - P^*}^2  
    \end{split}    
    \end{align}
    where the last inequality follows from the trace inequality in Lemma \ref{lm:traceIneq}.
\end{pf}

We are ready to state the main results on robustness analysis of the natural gradient descent algorithm in the presence of perturbations. The core idea is to demonstrate that the expression
\begin{align}
    \mathcal{V}_5 = \innprod{K-K^*}{K-K^*}_{Y^*} + \mathcal{J}_2(K) - \mathcal{J}_2(K^*)
\end{align}
serves as a small-disturbance ISS-Lyapunov function.

\begin{thm} \label{thm:ISSnature}
If $0 < \eta_{\min} \le \eta(k) \le \min\left\{\frac{1}{2\norm{R}},\frac{1}{6\norm{R}c(K(k))} \right\}$, the natural gradient descent algorithm in \eqref{eq:Gradient_NaturalCont} is small-disturbance ISS with respect to $W$.
\end{thm}

\begin{pf}
Define $\mathcal{V}_3(K)$ as
\begin{align}
    \mathcal{V}_3(K) = \innprod{K - K^*}{K - K^* }_{Y^*}.
\end{align}    
The difference of $\mathcal{V}_3(K)$ along the trajectories of \eqref{eq:Gradient_NaturalCont} is
\begin{align}
    &\mathcal{V}_3(K(k+1)) - \mathcal{V}_3(K(k)) \nonumber\\
    &= - 4\eta(k)\innprod{K(k) - K^* }{RK(k) - B^\top P(k)}_{Y^*} \nonumber\\
    &\quad + 4\eta(k)^2 \innprod{RK(k) - B^\top P(k)}{RK(k) - B^\top P(k) }_{Y^*}  \nonumber\\
    &\quad - 2\eta(k)\innprod{K(k) - K^*}{W(k)}_{Y^*} \nonumber\\
    &\quad + 4\eta(k)^2\innprod{RK(k) - B^\top P(k)}{W(k)}_{Y^*}  \nonumber\\
    &\quad +  \eta(k)^2\innprod{W(k)}{W(k)}_{Y^*} 
\end{align}
Lemma \ref{lm:PdiffNatural}, the trace inequality in Lemma \ref{lm:traceIneq}, the Cauchy-Schwarz inequality in Lemma \ref{lm:CSInequality}, and Young's inequality imply that
\begin{align}
    &\mathcal{V}_3(K(k+1)) - \mathcal{V}_3(K(k))  \nonumber\\
    &\le -2\eta(k) \Tr{P(k) - P^*}  \nonumber\\
    &\quad- \eta(k)\innprod{K(k) - K^* }{R(K(k) - K^*)}_{Y^*}  \nonumber\\
    &\quad + 6\eta(k)^2 \innprod{RK(k) - B^\top P(k)}{RK(k) - B^\top P(k) }_{Y^*}  \nonumber\\
    &\quad + \Big(3\eta(k) + \frac{1}{\eigmin{R}}\Big)\eta(k)\innprod{W(k)}{W(k)}_{Y^*} 
\end{align}
Considering Lemma \ref{lm:KK'bound} and applying the trace inequality in Lemma \ref{lm:traceIneq} again give that
\begin{align}
\begin{split}
    &\mathcal{V}_3(K(k+1)) - \mathcal{V}_3(K(k)) \\
    &\le -2\eta(k) \Tr{P(k) - P^*} \\
    &\quad - \eta(k)\innprod{K(k) - K^* }{R(K(k) - K^*)}_{Y^*} \\
    &\quad + 6\eta(k)^2 \norm{R} c(K(k)) \Tr{P(k) - P^*} \\
    &\quad + \Big(3\eta(k) + \frac{1}{\eigmin{R}}\Big)\eta(k)\innprod{W(k)}{W(k)}_{Y^*} \\
\end{split}
\end{align}
Hence, when $\eta(k) \le \frac{1}{6\norm{R}c(K(k))}$, it holds that
\begin{align}\label{eq:V1diff}
    &\mathcal{V}_3(K(k+1)) - \mathcal{V}_3(K(k)) \le - \eta(k) \Tr{P(k) - P^*} \nonumber \\
    &\quad - \eta(k)\innprod{K(k) - K^*}{R(K(k) - K^*)}_{Y^*} \nonumber\\
    &\quad  + \Big(3\eta(k) + \frac{1}{\eigmin{R}}\Big) \eta(k)\norm{Y^*}\norm{W(k)}_F^2.
\end{align}

We next derive the difference of $\mathcal{V}_4(K) = \mathcal{J}_2(K) - \mathcal{J}_2(K^*)$ along the trajectories of \eqref{eq:Gradient_NaturalCont}. At the $(k+1)$\textsuperscript{th} iteration, $P(k+1)$ is the solution of the Lyapunov equation
\begin{align}\label{eq:(k+1)thLyaunov}
\begin{split}
    &(A-BK(k+1))^\top  P(k+1) \\
    &\quad + P(k+1)(A-BK(k+1)) \\
    &\quad + Q + K(k+1)^\top R K(k+1) = 0.
\end{split}
\end{align} 
The Lyapunov equation at the $k$\textsuperscript{th} iteration can be rewritten as
\begin{align}\label{eq:kthLyaunov}
    &(A-BK(k+1))^\top  P(k) + P(k)(A-BK(k+1)) + Q \nonumber\\
    &\quad + K(k)^\top R K(k) + (K(k+1) - K(k))^\top B^\top  P(k) \nonumber\\
    &\quad  + P(k)B(K(k+1) - K(k)) = 0.
\end{align}
Subtracting \eqref{eq:kthLyaunov} from \eqref{eq:(k+1)thLyaunov} and considering \eqref{eq:Gradient_NaturalCont} give that
\begin{align}
\begin{split}
    0 &= (A-BK(k+1))^\top  (P(k+1) - P(k)) \\
    &+ (P(k+1) - P(k))(A-BK(k+1)) \\
    &+ (RK(k) - B^\top P(k))^\top (4\eta(k)^2R - 4\eta(k) I_m )\\
    &\quad (RK(k) - B^\top P(k)) + \eta(k)^2W(k)^\top  R W(k)\\
    & + (RK(k) - B^\top P(k))^\top (2\eta(k)^2 R - \eta(k)I_m )W(k) \\
    &+ W(k)^\top  (2\eta(k)^2 R - \eta(k)I_m)(RK(k) - B^\top P(k)) .
\end{split}
\end{align}
If $\eta(k) \le \frac{1}{2\norm{R}}$, $4\eta(k)^2 R - 4\eta(k) I_m \preceq 4\eta(k)^2\norm{R}I_m - 4\eta(k) I_m \preceq -2\eta(k) I_m$. By Lemma \ref{lm:traceIneq} and Corollary \ref{cor:traceLya}, it holds that
\begin{align}\label{eq:V2diff}
    &\mathcal{V}_4(K(k+1)) - \mathcal{V}_4(K(k)) \nonumber\\
    &\le -2\eta(k)\innprod{RK(k) - B^\top P(k)}{RK(k) - B^\top P(k)}_{Y(k+1)}  \nonumber\\
    &\quad+ \eta(k)^2 \innprod{W(k)}{RW(k)}_{Y(k+1)}   \nonumber\\
    &\quad+2\eta(k)\langle W(k), \nonumber\\
    &\qquad (2\eta(k) R - I_m )(RK(k) - B^\top P(k))\rangle_{Y(k+1)}  \nonumber\\
    &\le (1 + \norm{R}\eta(k))\eta(k) \frac{\mathcal{J}_2(K(k+1))}{\eigmin{Q}} \norm{W(k)}_F^2 
\end{align}
where Lemma \ref{lm:CSInequality}, Young's inequality, and \eqref{eq:traceYbound} are applied to derive the last inequality.

We are now prepared to demonstrate that the function
\begin{align}
    \mathcal{V}_5(K) = \mathcal{V}_3(K) + \mathcal{V}_4(K)
\end{align}
is indeed a small-disturbance ISS-Lyapunov function for the system in \eqref{eq:Gradient_NaturalCont}. Without loss of generality, assume that $\eigmin{R} \le 1$. By combining \eqref{eq:V1diff} and \eqref{eq:V2diff}, it follows that 
\begin{align}\label{eq:V3diff}
    &\mathcal{V}_5(K(k+1)) - \mathcal{V}_5(K(k)) \le - \eta(k) \eigmin{R} \mathcal{V}_5(K(k)) \nonumber\\
    &\quad + \eta(k)(c_1 + c_2 \mathcal{J}_2(K(k+1))) \norm{W(k)}_F^2
\end{align}
where $c_1 = \frac{3\eigmin{R} + 2 \norm{R}}{2\norm{R} \eigmin{R}}\norm{Y^*}$ and $c_2 = \frac{3}{2\eigmin{Q}}$. Since $\eta(k) \le \frac{1}{2\norm{R}}$ and $\mathcal{J}_2(K(k+1)) = \mathcal{V}_4(K(k+1))+ \mathcal{J}_2(K^*) \le \mathcal{V}_5(K(k+1))+ \mathcal{J}_2(K^*)$, it readily follows from \eqref{eq:V3diff} that
\begin{align}
    &\Big(1-\frac{\eta(k)\eigmin{R}}{2}\Big)(\mathcal{V}_5(K(k+1)) - \mathcal{V}_5(K(k))) \nonumber\\
    &\le - \frac{\eta(k)\eigmin{R}}{2}\mathcal{V}_5(K(k)) \nonumber\\
    &\quad - \frac{\eta(k)\eigmin{R}}{2}\mathcal{V}_5(K(k+1)) + \eta(k)(c_1+c_2\mathcal{J}_2(K^*)  \nonumber\\
    &\quad  + c_2 \mathcal{V}_5(K(k+1))) \norm{W(k)}_F^2
\end{align}
Hence, if 
\begin{align}
\begin{split}
    \norm{W(k)}_F^2 &\le \frac{\eigmin{R}}{2c_2} \frac{\mathcal{V}_5(K(k))}{1+\mathcal{V}_5(K(k))} \\
    &=: \sigma_1(\mathcal{V}_5(K(k))) < \frac{\eigmin{R}}{2c_2},     
\end{split}
\end{align}
it holds that
\begin{align}
\begin{split}
    &\Big(1-\frac{\eta(k)\eigmin{R}}{2}\Big)(\mathcal{V}_5(K(k+1)) - \mathcal{V}_5(K(k))) \\
    & \quad \le - \frac{\eta(k)\eigmin{R}}{2}\mathcal{V}_5(K(k)) \\
    &\qquad + \eta(k)(c_1+c_2\mathcal{J}_2(K^*)) \norm{W(k)}_F^2.
\end{split}
\end{align}
In addition, if 
\begin{align}
    \norm{W(k)}_F^2 \le  \frac{\eigmin{R}\mathcal{V}_5(K(k))}{4(c_1 + c_2 \mathcal{J}_2(K^*))} =: \sigma_2(\mathcal{V}_5(K(k))), 
\end{align}
it follows that
\begin{align}
\begin{split}
    &\mathcal{V}_5(K(k+1)) - \mathcal{V}_5(K(k)) \\
    &\quad \le -\frac{\eta(k) \eigmin{R}}{4-2\eta(k) \eigmin{R}} \mathcal{V}_5(K(k)) \\
    &\quad \le -\frac{\eta(k) \eigmin{R}}{4}\mathcal{V}_5(K(k))
\end{split}
\end{align}
In summary, if $\norm{W(k)}_F^2 \le \sigma(\mathcal{V}_5(K(k)))$, where $\sigma(r) = \min(\sigma_1(r), \sigma_2(r))$, $\mathcal{V}_5(K(k+1)) - \mathcal{V}_5(K(k)) \le -\frac{\eta(k) \eigmin{R}}{4}\mathcal{V}_5(K(k))$. Since $\sigma$ is a $\mathcal{K}$-function, according to Theorem \ref{thm:smallISSSufficient}, the system \eqref{eq:Gradient_NaturalCont} is small-disturbance ISS.
\end{pf}

\subsection{Small-Disturbance ISS of the Gauss-Newton Method}
This subsection analyzes the robustness of the Gauss-Newton method for solving the policy optimization problem of the LQR presented in \eqref{eq:costJc_closedform}. The action of the Hessian on $\de K \in \R^{m \times n}$ can be formulated as 
\begin{align}\label{eq:Hess2}
\begin{split}
    &\nabla^2 \mathcal{J}_2(K)[\de K] =\\
    & 2R\de KY_K  - 2B^\top \de P_KY_K + 2(RK - B^\top P_K)\de Y_K.
\end{split}
\end{align}
When $K = K^*$, we have $RK^* - B^\top P^* = 0$, and it follows from \eqref{eq:PKderivative} that $\de P_K = 0$ for all $\de K \in \R^{m \times n}$. Therefore, in the vicinity of $K^*$, the last two terms in \eqref{eq:Hess2} become negligible, allowing us to approximate the Hessian as
\begin{align}
     \nabla^2 \mathcal{J}_2(K)[\de K] \approx 2R\de KY_K.
\end{align}
This approximation of the Hessian is derived based on arguments similar to those used in the Gauss-Newton method \cite[Section 10.3]{book_Nocedal}. Hence, the update direction of the Gauss-Newton method is
$-(K-R^{-1}B^\top P_K)$, which is obtained by solving $\de K$ from
\begin{align}
    \nabla^2 \mathcal{J}_2(K)[\de K] \approx 2R\de KY_K = -\nabla \mathcal{J}_2(K).
\end{align}
Under the perturbation, the Gauss-Newton algorithm is
\begin{align}\label{eq:GradientFlow_Gauss-Newton}
    K(k+1) = K(k) -\eta(k)(K(k) - R^{-1}B^\top P(k) + W(k))
\end{align}

\begin{rem}
    The Gauss-Newton method in \eqref{eq:GradientFlow_Gauss-Newton} is derived based on the policy optimization of the LQR cost $\mathcal{J}_2(K)$. The update in \eqref{eq:GradientFlow_Gauss-Newton} can also be interpreted as an application of the classical Newton's method to solve the nonlinear algebraic Riccati equation (ARE) \cite{Kleinman1968,guo2000newton}:
\begin{align}
    \mathcal{R}(X) = A^\top X + XA - XBR^{-1}B^\top X + Q = 0. 
\end{align}    
\end{rem}

Indeed, the action of the gradient of the Riccati operator $\mathcal{R}(X)$ on $\de X \in \mathbb{S}^n$ is given by
\begin{align}
\begin{split}
    \nabla \mathcal{R}(X)[\de X] &= (A -BR^{-1}B^\top X)^\top \de X \\
    &+ \de X (A -BR^{-1}B^\top X).
\end{split}
\end{align}
According to Newton's method, at the $(k+1)$\textsuperscript{th} iteration, $X(k+1)$ is updated as
\begin{align}\label{eq:Newton'sUpdate}
    X(k+1) = X(k) + \eta(k)N(k)
\end{align}
where $N(k)$, the Newton updated direction, is the solution of $\nabla \mathcal{R}(X(k))[N(k)] = -\mathcal{R}(X(k))$, that is
\begin{align}\label{eq:NRiccati}
\begin{split}
    &(A -BR^{-1}B^\top X(k))^\top N(k) \\
    &\quad + N(k) (A -BR^{-1}B^\top X(k)) = -\mathcal{R}(X(k))    
\end{split}
\end{align}
Let us define $\bar{K}(k) = R^{-1}B^\top X(k)$ and $\bar{P}(k) = N(k) + X(k)$. With these definitions, \eqref{eq:NRiccati} can be reformulated as
\begin{align}
\begin{split}
    &(A -B\bar{K}(k))^\top \bar{P}(k) +\bar{P}(k) (A -B\bar{K}(k)) \\
    &\quad + Q + \bar{K}(k)^\top R\bar{K}(k)  = 0.
\end{split}
\end{align}
In addition, considering \eqref{eq:Newton'sUpdate} and the relation $N(k) = \bar{P}(k) - X(k)$, the recursive formula of $\bar{K}(k)$ becomes 
\begin{align}\label{eq:NewtonKupdate}
\begin{split}
    \bar{K}(k+1) &= R^{-1}B^\top X(k+1) \\
    &=\bar{K}(k) + \eta(k) R^{-1}B^\top(\bar{P}(k) - X(k)) \\
    &=\bar{K}(k) - \eta(k)(\bar{K}(k) - R^{-1}B^\top \bar{P}(k)). 
\end{split}
\end{align}
Hence, \eqref{eq:NewtonKupdate} is equivalent to \eqref{eq:GradientFlow_Gauss-Newton} without perturbation, which implies that the Gauss-Newton method in \eqref{eq:GradientFlow_Gauss-Newton} coincides with Newton's method in \eqref{eq:Newton'sUpdate} for solving the ARE. This interpretation establishes a connection between the Gauss-Newton method for policy optimization and the classical Newton's method for solving the ARE. The following theorem shows the small-disturbance ISS property of the Gauss-Newton method in \eqref{eq:GradientFlow_Gauss-Newton}.

\begin{thm}  \label{thm:smallISSGauss-Newton}
If $0 < \eta_{\min} \le \eta(k)  \le \min\left\{1, \frac{1}{4c(K(k))}\right\}$, the Gauss-Newton method in \eqref{eq:GradientFlow_Gauss-Newton} is small-disturbance ISS with respect to $W$.
\end{thm}
\begin{pf}
It can be shown that $\mathcal{V}_6(K) = \mathcal{V}_4(K) + \frac{1}{2}\innprod{K - K^*}{R(K - K^*)}_{Y^*}$ is a small-disturbance ISS-Lyapunov function by following the proof of Theorem \ref{thm:ISSnature}.
\end{pf}

\section{Simulations}
\begin{figure}[tb]
    \begin{subfigure}{.49\linewidth}
    \centering
      	\includegraphics [width=0.98\linewidth]{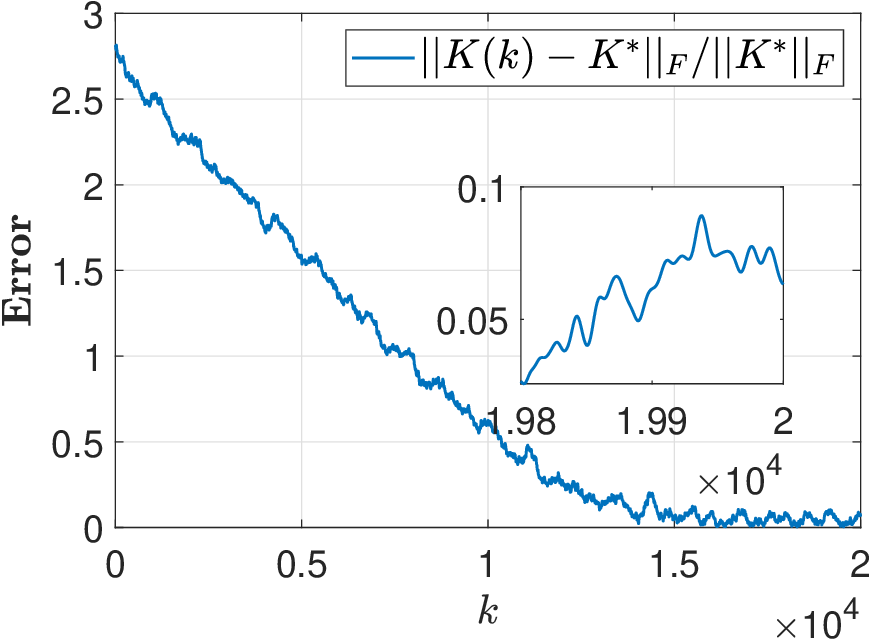}
    	\label{fig: GradGain}  
    \end{subfigure}
    \begin{subfigure}{.49\linewidth}
    	\centering
    	\includegraphics [width=0.98\linewidth]{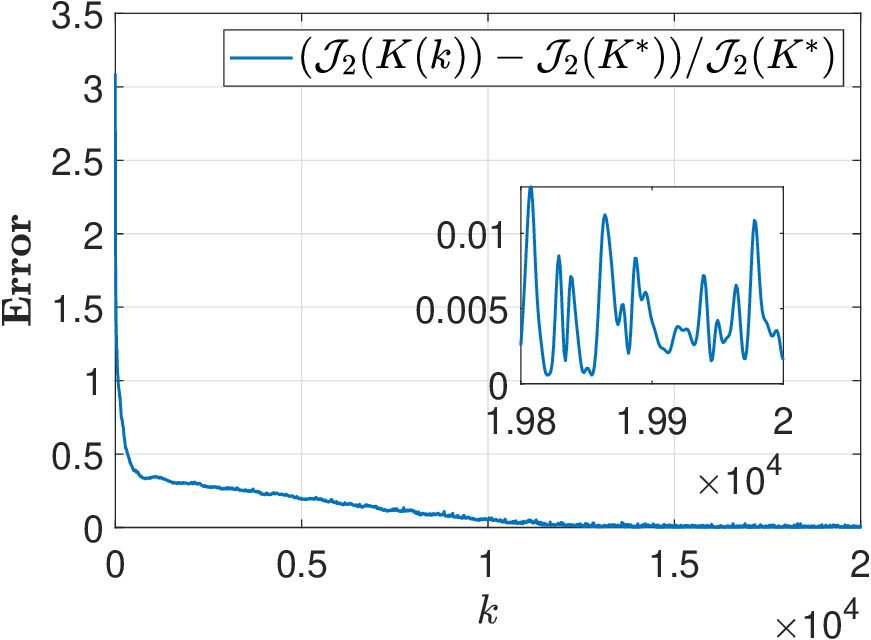}
    	\label{fig: GradCost}
    \end{subfigure}	
	\caption{Robustness test of the perturbed standard gradient descent \eqref{eq:graddescentLQR} when $\sigma = 1$.}
	\label{fig:GradFlow}
\end{figure}

\begin{figure}[tb]
    \begin{subfigure}{.49\linewidth}
    \centering
      	\includegraphics [width=0.98\linewidth]{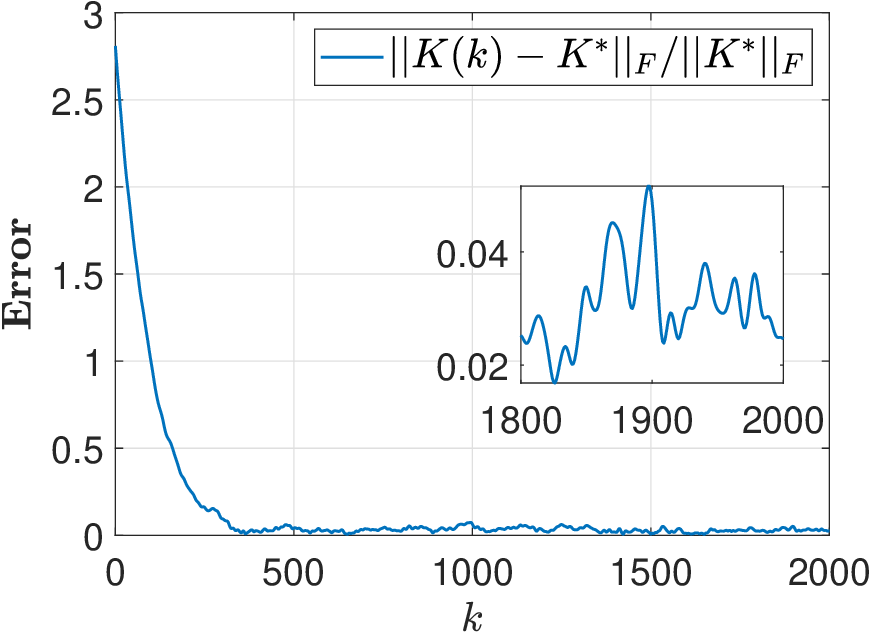}
    	\label{fig: NaturalGradGain}  
    \end{subfigure}
    \begin{subfigure}{.49\linewidth}
    	\centering
    	\includegraphics [width=0.98\linewidth]{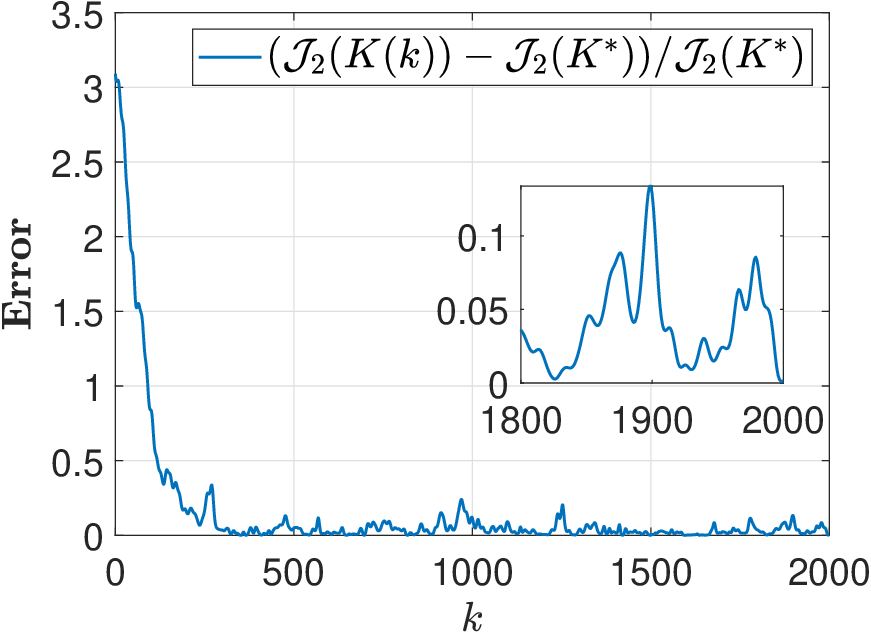}
    	\label{fig: NaturalGradCost}
    \end{subfigure}	
	\caption{Robustness test of the perturbed natural gradient descent \eqref{eq:Gradient_NaturalCont} when $\sigma = 1$.}
	\label{fig:NaturalGradFlow}
\end{figure}

\begin{figure}[tb]
    \begin{subfigure}{.49\linewidth}
    \centering
      	\includegraphics [width=0.98\linewidth]{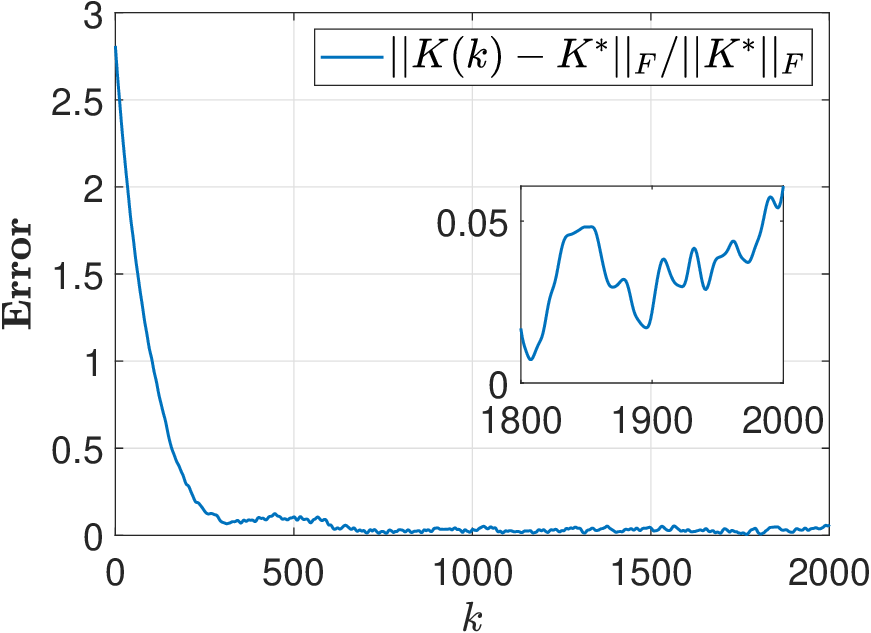}
    	\label{fig: NewtonGradGain}  
    \end{subfigure}
    \begin{subfigure}{.49\linewidth}
    	\centering
    	\includegraphics [width=0.98\linewidth]{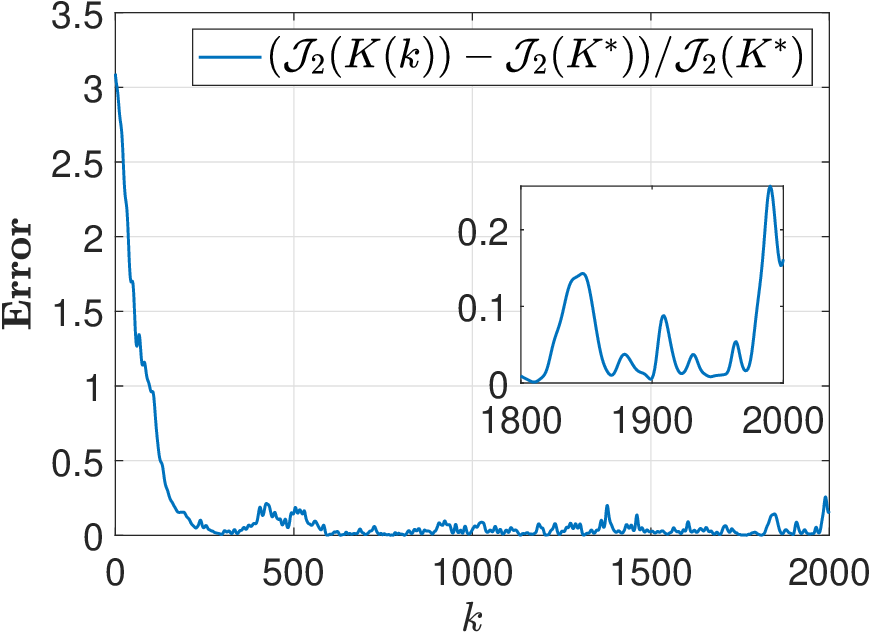}
    	\label{fig: NewtonGradCost}
    \end{subfigure}	
	\caption{Robustness test of the perturbed Gauss-Newton gradient descent \eqref{eq:GradientFlow_Gauss-Newton} when $\sigma = 1$.}
	\label{fig:NewtonGradFlow}
\end{figure}

\Review{
In this section, the theoretical results are illustrated by an example known as policy optimization for a cart-pole system. The mass of the cart is $m_c = 1$~kg, the mass of the pendulum is $m_p = 0.1$~kg, the length between the center of mass of the pendulum and the cart is $l_p = 0.5$~m, and the gravitational acceleration is $g_c = 9.8$~m/sec\textsuperscript{2}. The system is linearized around the equilibrium and its state-space representation is
\begin{align}
    \dot{x} = \begin{bmatrix} 0	 &1 &0 &0 \\
    0 &0 &-0.717 &0	\\
    0 &0 &0	&1 \\
    0 &0 &15.77	&0
\end{bmatrix}x + \begin{bmatrix}
    0 \\
    0.98 \\
    0 \\
    -1.46
\end{bmatrix} u.
\end{align}
The weighting matrices of the objective function \eqref{eq:LQRcost} are $Q=10I_4$ and $R=1$. The initial control gain is 
\begin{align}
    K(0) = \begin{bmatrix}-58.6 &-43.6 &-167.8 &-43.5 \end{bmatrix}.
\end{align}
Each element of $W(k)$ is set as
\begin{align}
    [W(k)]_{i,j} = \sigma\sum_{l=1}^{100}\sin( \omega_{l}^{i,j} k),
\end{align}
where $[W(k)]_{i,j}$ is the element of $W(k)$ at the $i$th row and $j$th column, and $\omega_{l}^{i,j}$ is independently and randomly sampled from the uniform distribution over $[-50,50]$. The step size is set as $\eta(k)=0.01$. The standard gradient descent algorithm terminates after $2\times10^4$ iterations or when the updated feedback gain fails to stabilize the closed-loop system. The natural gradient descent and Gauss–Newton methods terminate after $2\times10^3$ iterations or when the updated feedback gain is not stabilizing.
\\
\\
When the noise level is small ($\sigma=1$), the trajectory of the perturbed standard gradient descent is shown in Fig. \ref{fig:GradFlow}. It is seen that under the perturbation $W$, the control gain approaches the optimal solution at around $k=1.5\times 10^4$, and finally stays within a neighborhood of the optimal solution $K^*$.  The trajectories of the perturbed natural gradient descent and Newton gradient descent are shown in Fig. \ref{fig:NaturalGradFlow} and Fig. \ref{fig:NewtonGradFlow}, respectively. It is shown that after around $k=500$, the trajectories of the perturbed natural gradient descent and Newton gradient descent approach the optimum and remain close to the optimal solution. 
}

\Review{
When the noise level is large ($\sigma = 20$), all three algorithms terminate before reaching the maximum number of iterations, indicating that the updated control gains are driven outside the admissible set. Consequently, the gradient-based algorithms are no longer able to stabilize the closed-loop system. These simulation results are consistent with the small-disturbance ISS properties established in Theorems \ref{thm:LQRgradientISS}, \ref{thm:ISSnature}, and \ref{thm:smallISSGauss-Newton}.
}

\begin{figure}[tb]
    \begin{subfigure}{.49\linewidth}
    \centering
      	\includegraphics [width=0.98\linewidth]{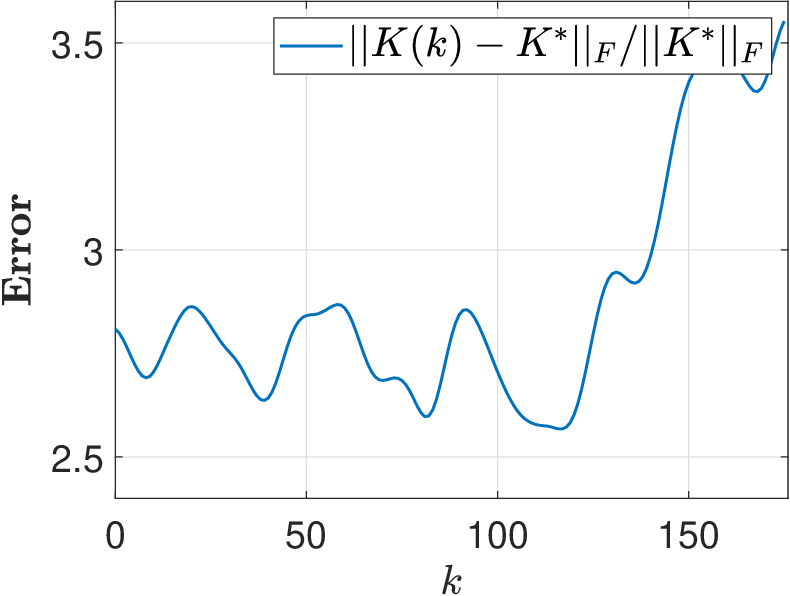}
    	\label{fig: GradGainLargeNoise}  
    \end{subfigure}
    \begin{subfigure}{.49\linewidth}
    	\centering
    	\includegraphics [width=0.98\linewidth]{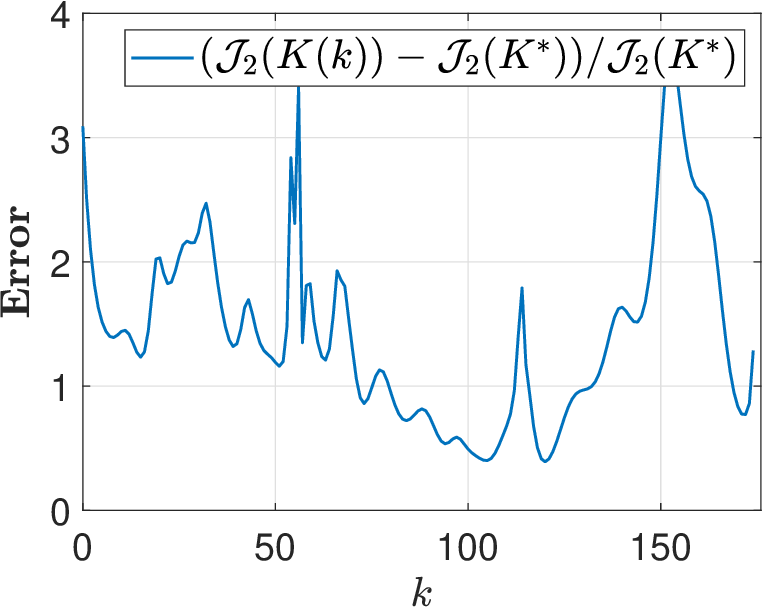}
    	\label{fig: GradCostLargeNoise}
    \end{subfigure}	
    \caption{Robustness test of the perturbed standard gradient descent algorithm \eqref{eq:graddescentLQR} with noise level $\sigma = 20$. The iteration terminates when the updated feedback gain fails to stabilize the closed-loop system.}
	\label{fig:GradFlowLargeNoise}
\end{figure}

\begin{figure}[tb]
    \begin{subfigure}{.49\linewidth}
    \centering
      	\includegraphics [width=0.98\linewidth]{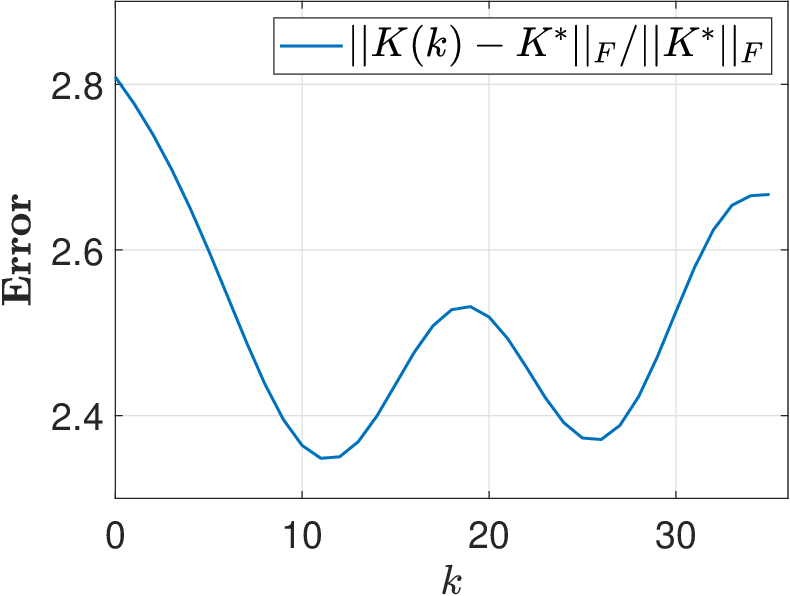}
    	\label{fig: NaturalGradGainLargeNoise}  
    \end{subfigure}
    \begin{subfigure}{.49\linewidth}
    	\centering
    	\includegraphics [width=0.98\linewidth]{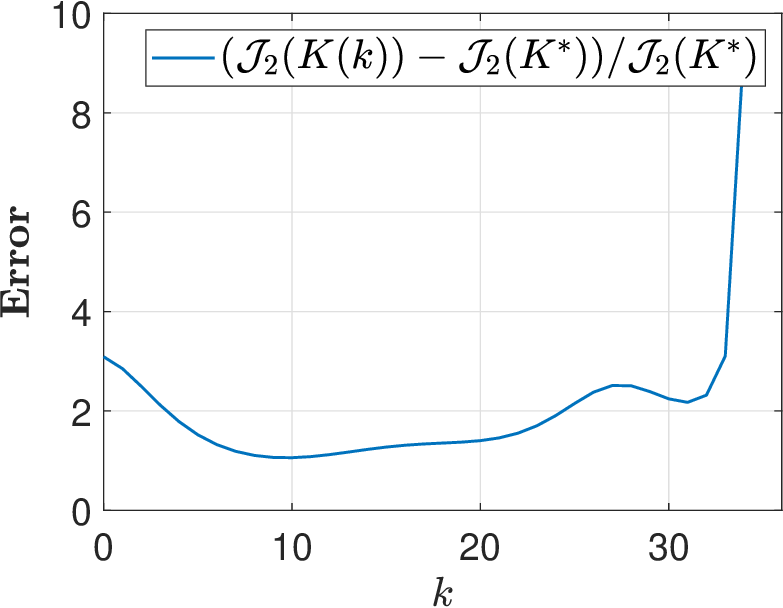}
    	\label{fig: NaturalGradCostLargeNoise}
    \end{subfigure}	
	\caption{Robustness test of the perturbed natural gradient descent algorithm \eqref{eq:Gradient_NaturalCont} with noise level $\sigma = 20$. The iteration terminates when the updated feedback gain fails to stabilize the closed-loop system.}
	\label{fig:NaturalGradFlowLargeNoise}
\end{figure}

\begin{figure}[tb]
    \begin{subfigure}{.49\linewidth}
    \centering
      	\includegraphics [width=0.98\linewidth]{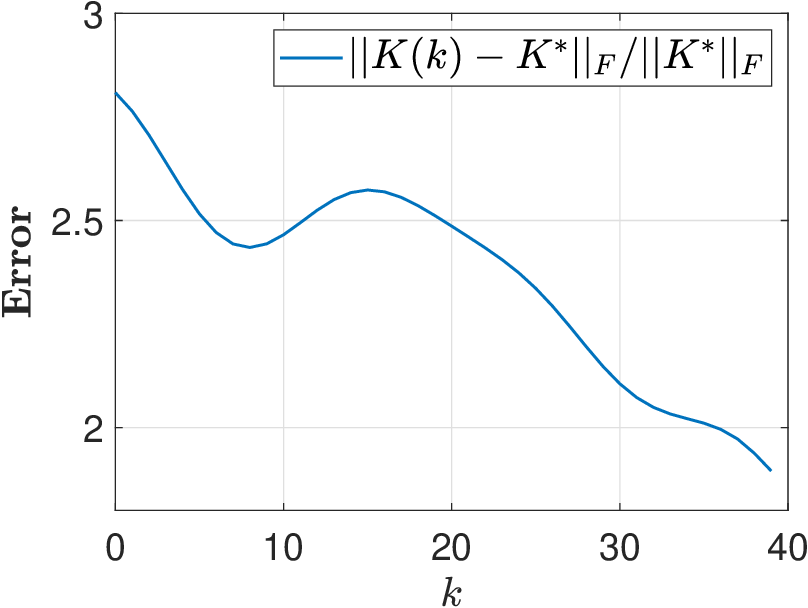}
    	\label{fig: NewtonGradGainLargeNoise}  
    \end{subfigure}
    \begin{subfigure}{.49\linewidth}
    	\centering
    	\includegraphics [width=0.98\linewidth]{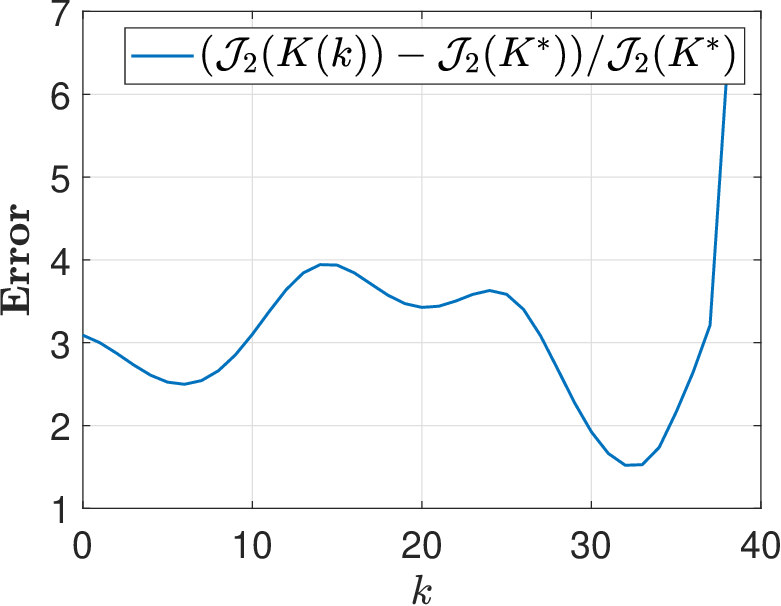}
    	\label{fig: NewtonGradCostLargeNoise}
    \end{subfigure}	
	\caption{ Robustness test of the perturbed Gauss-Newton gradient descent algorithm \eqref{eq:GradientFlow_Gauss-Newton} with noise level $\sigma = 20$. The iteration terminates when the updated feedback gain fails to stabilize the closed-loop system.}
	\label{fig:NewtonGradFlowLargeNoise}
\end{figure}

\section{Conclusions}
This article introduces the concept of small-disturbance ISS as a unified framework for analyzing the robustness of gradient descent algorithms. Small-disturbance ISS provides a systematic approach to quantify the transient behavior, convergence speed, and robustness of gradient descent algorithms under perturbations. By generalizing the classical linear PL condition to a nonlinear version, referred to as the $\mathcal{K}$-PL condition, we show that gradient descent algorithms are small-disturbance ISS, provided the objective function satisfies the $\mathcal{K}$-PL condition. As a direct application to LQR, we demonstrate that three popular policy gradient algorithms in RL--standard policy gradient, natural policy gradient, and Gauss-Newton method--are all small-disturbance ISS.

\Review{
Future work will extend the small-disturbance ISS framework to broader classes of gradient-based optimization algorithms, including constrained optimization problems, projected gradient methods, and proximal gradient algorithms for nonsmooth objectives. Another important direction is to study ISS properties of policy optimization algorithms in reinforcement learning for general Markov decision processes. The proposed framework may provide a systematic tool to quantify robustness margins and guide the design of reliable reinforcement learning algorithms.
}
\bibliographystyle{plain}        
\bibliography{reference}           %

\end{document}